\documentclass{cas-sc}
\pdfoutput=1

\usepackage{appendix}
\usepackage{graphbox}
\usepackage{mathtools}
\usepackage{dashrule}
\usepackage{enumitem}
\usepackage[numbers]{natbib}

\definecolor{stability_dark}{RGB}{102, 102, 102}
\definecolor{stability_light}{RGB}{204, 204, 204}

\definecolor{matlab_blue}{RGB}{0, 114, 189}
\definecolor{matlab_red}{RGB}{217, 83, 25}
\definecolor{matlab_yellow}{RGB}{237, 177, 32}
\definecolor{matlab_purple}{RGB}{126, 47, 142}
\definecolor{matlab_green}{RGB}{119, 172, 48}
\definecolor{matlab_azure}{RGB}{77, 190, 238}
\definecolor{matlab_dark_red}{RGB}{162, 20, 47}

\setcitestyle{square}
\setcitestyle{citesep={,}}
\hypersetup{
 colorlinks=true,
 citecolor=Cerulean,
 linkcolor=Cerulean,
 urlcolor=Cerulean}

\begin{document}

\title[mode = title]{On the Stability of Exponential Integrators for Non-Diffusive Equations}
\tnotemark[1]
\tnotetext[1]{This work was funded by the National Science Foundation, Computational Mathematics Program DMS-2012875}

\shorttitle{On the Stability of Exponential Integrators for Non-Diffusive Equations}
\shortauthors{T Buvoli et~al.}

\author[1]{Tommaso Buvoli}
\cormark[1]
\ead{tbuvoli@ucmerced.edu}
\address[1]{University of California, Merced, Merced CA 95343, USA}

\author[2]{Michael Minion}
\ead{mlminion@lbl.gov}
\address[2]{Lawrence Berkeley National Lab,Berkeley, CA 94720, USA}

\cortext[cor1]{Corresponding author}

\begin{abstract}
	Exponential integrators are a well-known class of time integration methods that have been the subject of many studies and developments in the past two decades. Surprisingly, there have been limited efforts to analyze their stability and efficiency on non-diffusive equations to date. In this paper we apply linear stability analysis to showcase the poor stability properties of exponential integrators on non-diffusive problems. We then propose a simple repartitioning approach that stabilizes the integrators and enables the  efficient solution of stiff, non-diffusive equations. To validate the effectiveness of our approach, we perform several numerical experiments that compare partitioned exponential integrators to unmodified ones. We also compare repartitioning to the well-known approach of adding hyperviscosity to the equation right-hand-side. Overall, we find that the repartitioning restores convergence at large timesteps and, unlike hyperviscosity, it does not require the use of high-order spatial derivatives.

\end{abstract}

\maketitle

\section{Introduction}
\label{sec:introduction}

In the past two decades, exponential integrators \cite{hochbruck2010exponentialreview} have emerged as an attractive alternative to both fully implicit and semi-implicit methods \cite{grooms2011IMEXETDCOMP, KassamTrefethen05ETDRK4, loffeld2013comparative, montanelli2016solving}. The growing interest in exponential integrators is jointly driven by new method families \cite{beylkin1998ELP, cox2002ETDRK4, KassamTrefethen05ETDRK4, krogstad2005IF,hochbruckostermann2005ETDRKSTIFFA, hochbruckostermann2005ETDRKSTIFFB, koikari2005rooted, hochbruck2010exponentialreview, ostermann2006general,buvoli2019esdc,buvoli2021epbm} and by advances in algorithms for computing the associated exponential functions \cite{al2010new,al2011computing,caliari2021accurate,haut2015,caliari2016leja,hochbruck1997krylov,NiesenWright2012Krylov,NiesenWright2011Krylov,GAUDREAULT2018236,caliari2014comparison,higham2020catalogue}. 

The purpose of this paper is to study the stability limitations of exponential integrators on problems with no diffusion, and to propose a repartitioning strategy that improves stability without damaging accuracy. At first glance, investigating the stability of exponential integrators may seem unfruitful since the methods have already been used to successfully solve a range of stiff partial differential equations including purely dispersive ones \cite{KassamTrefethen05ETDRK4,montanelli2016solving,grooms2011IMEXETDCOMP}. However, as was recently described in \cite{crouseilles2020exponential}, exponential integrators possess linear stability regions for non-diffusive problems that are on par with those of explicit methods. These instabilities are often very small in magnitude, which explains why they can sometimes go unnoticed. %
	In this work we will present two simple equations that excite instabilities in exponential integrators, and then demonstrate how a simple repartitioning approach can be used to stabilize them.

Our motivation for pursuing this work goes beyond simply stabilizing integrators. In particular, we believe that incorporating exponential integrators into parallel-in-time frameworks has the potential to produce new and efficient parallel methods for solving large scale application problems. These methods will be of particular interest if exponential integrators can also be used to achieve meaningful parallel speedup for equations with little to no diffusion.  Unfortunately, we have found that the effectiveness of exponential integrators in parallel frameworks is significantly hindered by the stability issues discussed in this paper. Therefore, in this preliminary work, we propose a strategy for overcoming this problem.

In light of this goal we will focus on three integrator families that can be immediately used within parallel frameworks. In particular, we will analyze: (1) exponential Runge-Kutta (ERK) methods \cite{krogstad2005generalized,cox2002ETDRK4}, which can be used to construct exponential Parareal methods, (2) exponential spectral deferred correction methods (ESDC) \cite{buvoli2019esdc}, which are the exponential equivalents of the classical SDC methods used in the Parallel Full Approximation Scheme in Space and Time (PFASST) framework \cite{emmett2012toward}, and (3) exponential polynomial block methods (EPBMs) \cite{buvoli2021epbm} that allow for parallel function evaluations. %
	Lastly, we note that the stability issues we describe affect many other families of exponential integrators, and therefore the developments from this paper can also be used to improve the efficiency of any serial exponential integrators for hyperbolic and non-diffusive problems.

The outline of this paper is as follows. In Section \ref{sec:exponential_integrators} we provide a brief overview of exponential integrators and the three families of methods that will be the subject of this work. In Section \ref{sec:model_problem} we then introduce a simple one-dimensional, nonlinear partial differential equation that leads to severe stepsize restrictions for exponential integrators. Next, in Section \ref{sec:linear_stability}, we use linear stability analysis to showcase the poor stability properties of exponential integrators on non-diffusive equations. Lastly, in Section \ref{sec:repartitioning}, we propose a simple repartitioning approach to stabilize exponential integrators, and compare it against the commonly applied practice of adding hyperviscosity.

\section{Exponential integrators}
\label{sec:exponential_integrators}

Exponential integrators \cite{hochbruck2010exponentialreview} are a family of time integration methods with coefficients that can be written in terms of exponential functions of a linear operator. They have proven to be highly suitable for solving stiff systems and, in certain situations, they can outperform both fully-implicit and linearly-implicit methods %
\cite{grooms2011IMEXETDCOMP, KassamTrefethen05ETDRK4, loffeld2013comparative, montanelli2016solving,hochbruck1997krylov}. 
Nearly all exponential integrators can be classified into two subfamilies: unpartitioned and partitioned. Unpartitioned exponential integrators \cite{hochbruck2009exponential,tokman2006efficient,tokman2011new,rainwater2014new,buvoli2019esdc,buvoli2021epbm} can be used to solve the generic differential equation $\mathbf{y}'=F(\mathbf{y})$ and require the exponentiation of the local Jacobian of $F(y)$ at each timestep. In contrast, partitioned exponential integrators \cite{beylkin1998ELP,cox2002ETDRK4,krogstad2005IF,hochbruckostermann2005ETDRKSTIFFB,buvoli2019esdc,buvoli2021epbm} are designed to solve the semilinear initial value problem
\begin{align}
	\mathbf{y}' = \mathbf{Ly} + N(t,\mathbf{y}), \quad \mathbf{y}(t_0) = \mathbf{y}_0,
	\label{eq:semilinear_ode}
\end{align}
and only require exponential functions involving the autonomous linear operator $\mathbf{L}$. 

In this work we will focus exclusively on partitioned exponential integrators. All partitioned exponential integrators can be derived by applying the variation of constants formula to (\ref{eq:semilinear_ode}), yielding
\begin{align}
  \mathbf{y}(t_0 + \Delta t) = e^{\Delta t\mathbf{L}}\mathbf{y}_0 + \int^{t_0+\Delta t}_{t_0}e^{(t-\tau)\mathbf{L}}
  N(\tau,\mathbf{y}(\tau))d\tau,
  \label{eq:exp_derivation_formula}
\end{align}
and then approximating the nonlinear term using a polynomial. For example, to obtain the first-order partitioned exponential Euler method, 
\begin{align}
	\mathbf{y}_{n+1} &= e^{\Delta t\mathbf{L}} \mathbf{y}_n + \int^{1}_0 e^{(1-s)\Delta t\mathbf{L}} N(t_n, \mathbf{y}_n) ds,
\end{align}
with stepsize $\Delta t$, we approximate the nonlinear term $N(s,\mathbf{y})$ in (\ref{eq:exp_derivation_formula}) with a zeroth-order polynomial $p(s) = N(t_n, y_n)$ and then rewrite the integrand in local coordinates using ${\tau = \Delta t s + t_n}$. Higher order exponential multistep, Runge-Kutta, and general linear methods can be constructed by approximating $N$ with a higher-order polynomial that respectively uses previous solution values, previous stage values, or any combination of the two.

In this work we will be concerned with three specific families of exponential integrators.
\begin{enumerate}
	\item {\em Exponential Runge-Kutta (ERK)} methods from \cite{cox2002ETDRK4,krogstad2005IF} accept a single input, compute $s$ stage values, and produce a single output. These methods attempt to achieve a high order of accuracy in the fewest number of stages, and the coefficients have been derived by satisfying nonlinear order conditions. In this work, we will consider the fourth-order exponential Runge-Kutta method ETDRK4-B from \cite{krogstad2005IF} which will be referred to simply as ERK4. %
	\item {\em Exponential spectral deferred correction (ESDC)} methods \cite{buvoli2019esdc} are a class of arbitrary order time integration schemes that iteratively improve a provisional solution by solving an integral equation that governs the error. All ESDC methods can be written as ERK methods with a large number of stages. In this paper we will consider a 6th-order ESDC method with four Gauss-Lobatto nodes that takes six correction iterations; we will refer to this method as ESDC6.	 
	\item {\em Exponential polynomial block methods (EPBM)} \cite{buvoli2021epbm} are multivalued, exponential general linear methods that advance a set of $q$ different solution values at each timestep. The input values are approximations to the solution at different time points, and the outputs are computed by approximating the nonlinearity in (\ref{eq:semilinear_ode}) using a high-order polynomial that interpolates the nonlinearity at the input values. The methods can be constructed at any order of accuracy and allow for parallel function evaluations.	Here we consider a fifth-order composite EPBM that is constructed using Legendre nodes and is run using an extrapolation factor of $\alpha = 1$.
\end{enumerate}

The formulas for all three exponential integrators are contained in Appendix \ref{app:coefficients}, and  Matlab code for running all the numerical experiments can be downloaded from \cite{BuvoliExponentialStabilityGithub}.
\section{A motivating example}
\label{sec:model_problem}

We begin by showcasing a simple one dimensional nonlinear partial differential equation that causes stability problems for exponential integrators. Specifically, we consider the zero-dispersion Schr\"{o}dinger equation (ZDS) 
\begin{align}
	\begin{aligned}
		& iu_t + i u_{xxx} + 2u|u|^2 = 0 \\
		& u(x, t=0) = 1 + \tfrac{1}{100} \exp(3ix / 4).
	\end{aligned}
	\label{eq:zds}
\end{align}
This complex-valued, dispersive partial differential equation models optical pulses in zero dispersion fibers \cite{agrawal2000nonlinear} and is obtained by replacing the second derivative in the canonical nonlinear Schr\"{o}dinger equation with an $iu_{xxx}$ term.

We equip the ZDS equation with periodic boundary conditions in $x$ and integrate to time $t=40$ using a Fourier pseudo-spectral method. To simplify the computation of the exponential functions, we solve the equation in Fourier space where the linear derivative operators are diagonal and (\ref{eq:zds}) reduces to the semilinear equation (\ref{eq:semilinear_ode}) with
	\begin{align*}
		\mathbf{L} = \text{diag}(i \mathbf{k}^3) 
		\quad \text{and} \quad
		N(\hat{u}) = 2 \mathcal{F}(\text{abs}(\mathcal{F}^{-1}(\hat{u})) .* \mathcal{F}^{-1}(\hat{u}))
	\end{align*}
	where $\mathbf{k}$ is a vector of Fourier wavenumbers and $\mathcal{F}$ denotes the discrete Fourier transform. The full numerical parameters are summarized in Table \ref{tab:zds_numerical_parameters}.

\begin{table}[h!]
	\centering
	\renewcommand*{\arraystretch}{1.5}
		\begin{tabular}{rl}
		{\em Domain:} & $x\in[-4\pi, 4\pi]$, $t \in [0, 40]$\\
		{\em Boundary conditions:} & periodic \\
		{\em Spatial discretization:} & Fourier pseudo-spectral \\
		{\em Number of spatial grid points:} & $N_x = 128$ \\
		{\em Dealiasing:} &  $\tfrac{3}{2}$ rule (see for example \cite[p. 84]{canuto2012spectral})
	\end{tabular}	
	\caption{Numerical parameters for the ZDS equation (\ref{eq:zds}).}
	\label{tab:zds_numerical_parameters}
	
\end{table}

\begin{figure*}[h!]
	
	\begin{tabular}{cc}
		ERK4 -- Physical & IMRK4 -- Physical \\
		\includegraphics[width=0.45\linewidth,trim={0 22 13 50},clip]{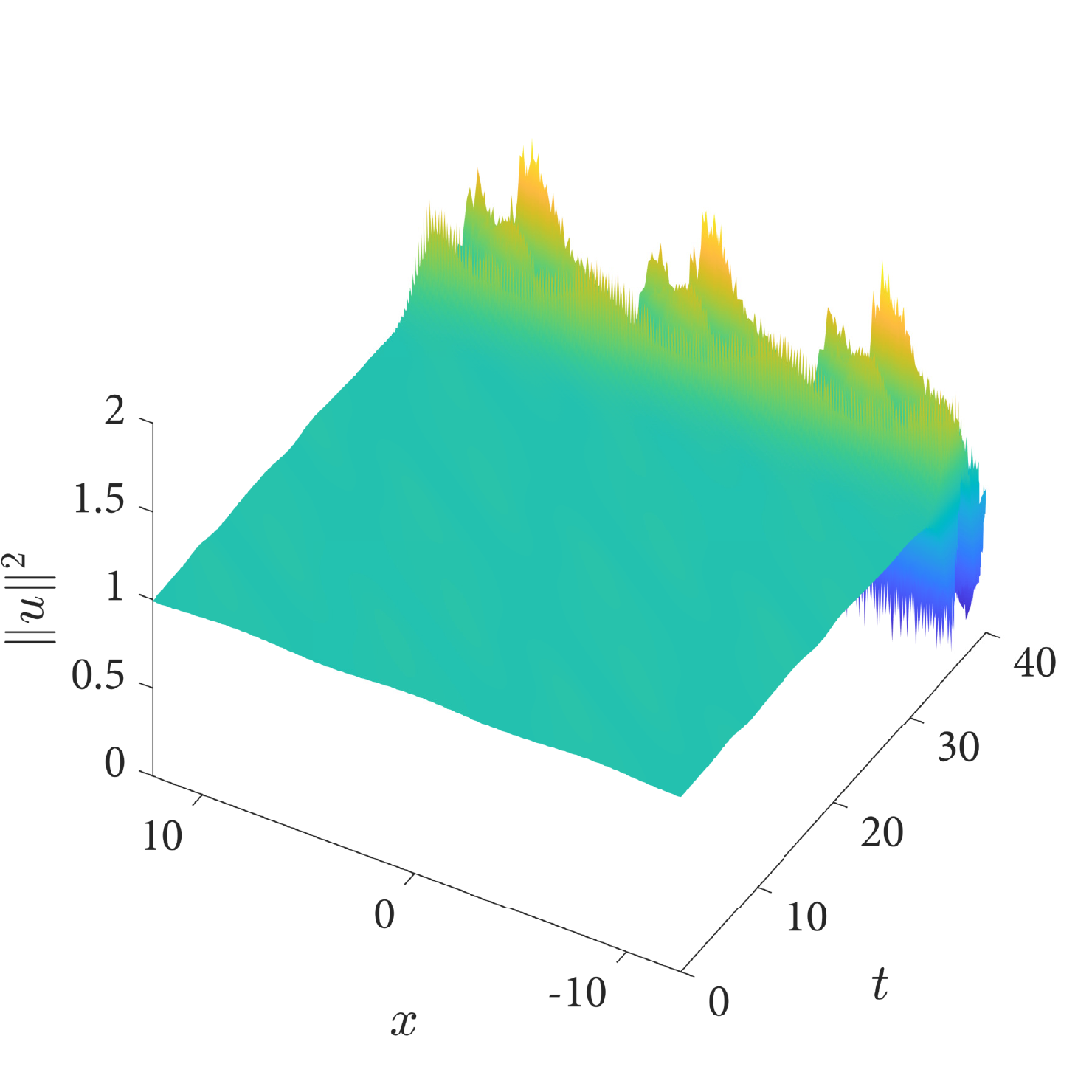} & 
		\includegraphics[width=0.45\linewidth,trim={0 22 13 50},clip]{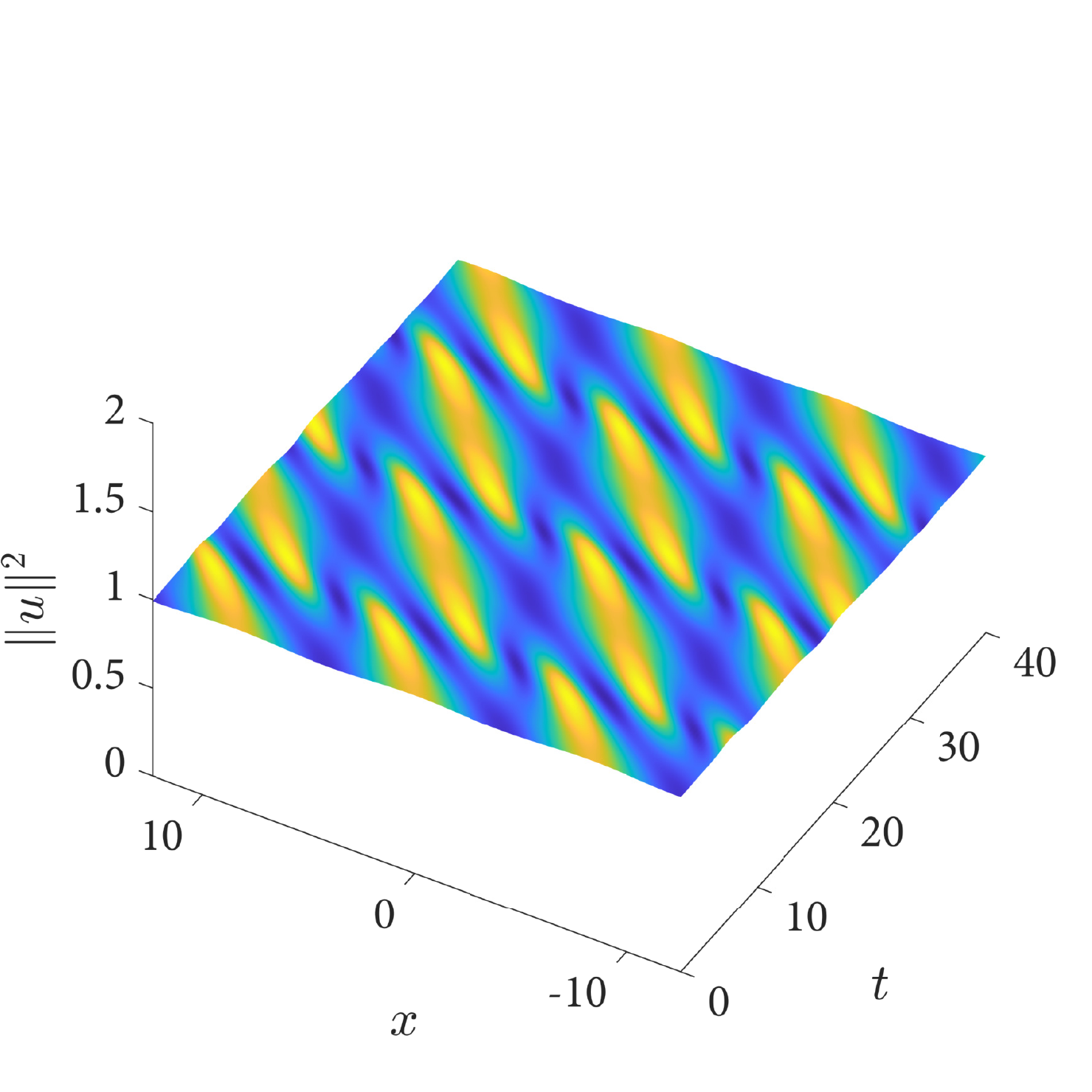} \\[1em]
		ERK4 -- Fourier & IMRK4 -- Fourier	 \\
		\includegraphics[width=0.45\linewidth,trim={0 22 13 50},clip]{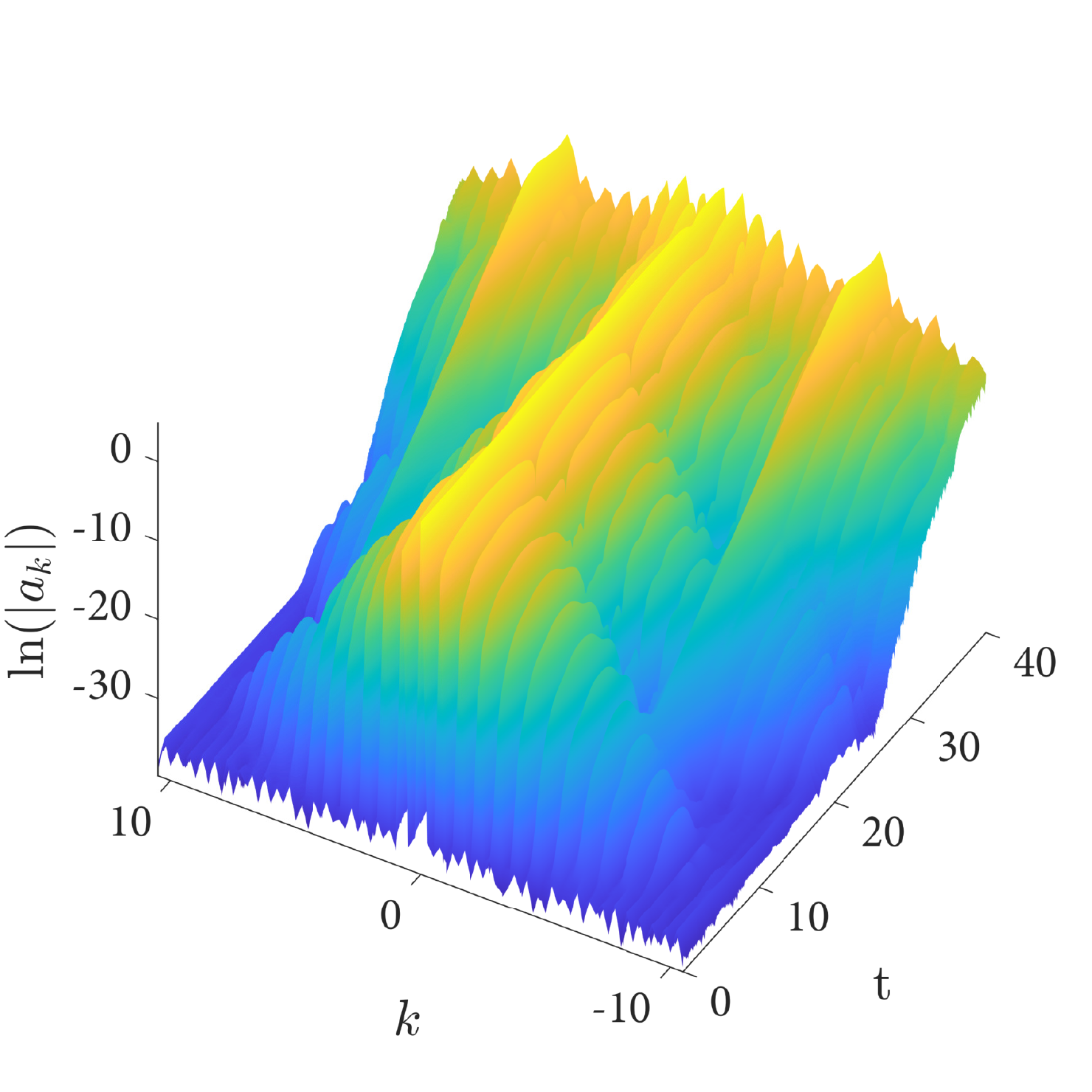} &
		\includegraphics[width=0.45\linewidth,trim={0 22 13 50},clip]{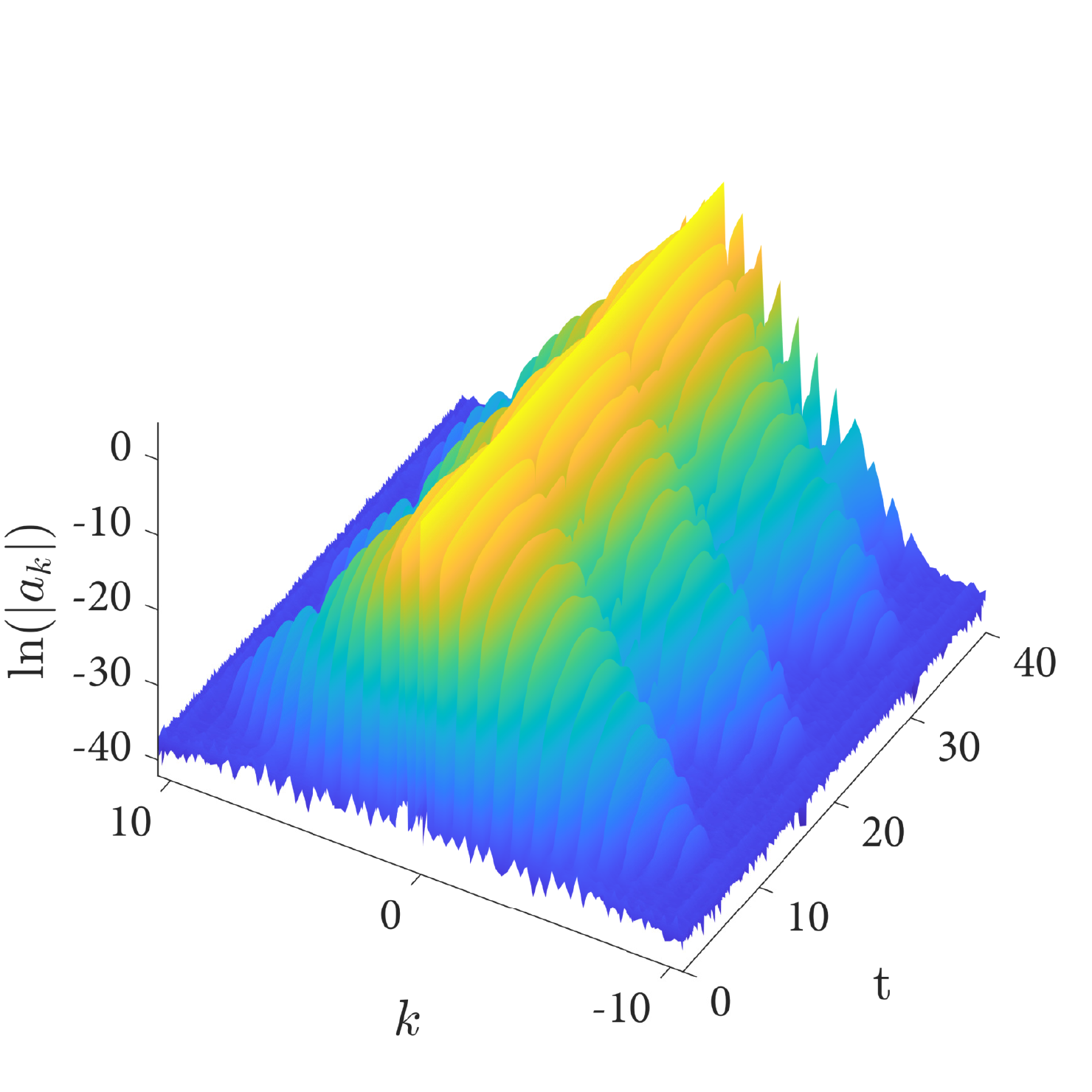}
	\end{tabular}
	\caption{Solution to the ZDS equation in physical and Fourier space using the ERK4 and the IMRK4 integrators that are both run using 2000 timesteps. In physical space the square of the two-norm of the solution is plotted, while in Fourier space we show the log of the absolute value of the Fourier coefficients.}
	\label{fig:zds-convergence}
	
\end{figure*}

In Figure \ref{fig:zds-convergence} we show the solution obtained by integrating in time using $2000$ steps of the ERK4 method. For comparison we also compute the solution using the fourth-order implicit-explicit method from \cite{kennedy2003additive} named ARK4(3)6L[2]SA, which we will simply refer to as IMRK4. Instabilities are clearly visible in the ERK4 solution while the IMRK4 solution is clean. Even more surprisingly, the problem is only mildly stiff. By using 2000 timesteps, or equivalently ${h=2\times 10^{-2}}$, the spectral radius of the linear operator $\rho(h\mathbf{L})$ is only $23.1525$; note that we are computing the spectral radius for the lower 2/3 of modes where antialiasing is not applied.

\begin{figure}[h!]
	\centering
	{\bf ZDS Convergence Diagram}
	
	\includegraphics[width=0.55\linewidth]{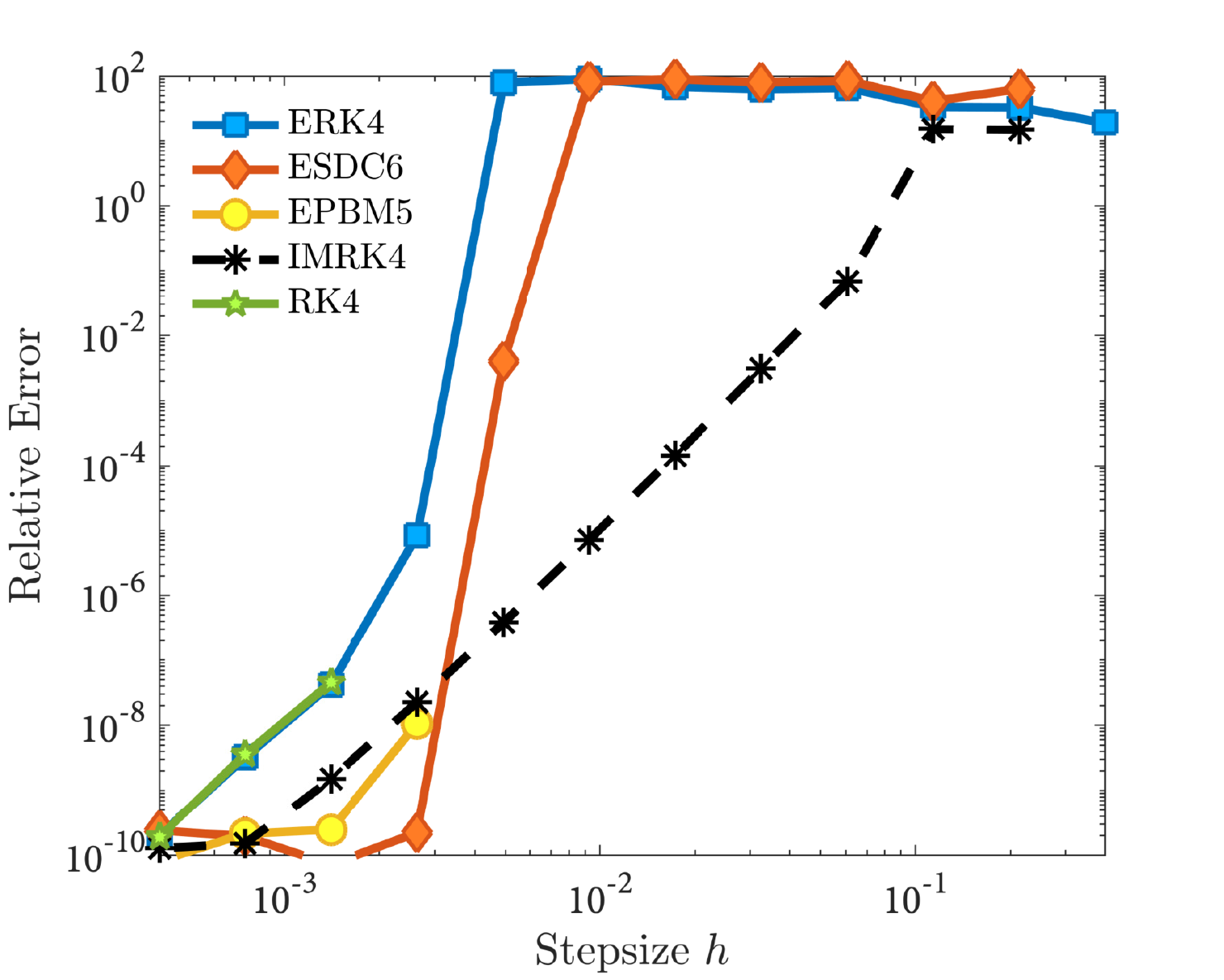}
	\caption{A convergence diagram for ERK4, ESDC6, and EPBM5 methods on the ZDS equation. We also include the IMRK4 method and classical explicit fourth-order RK4. At coarse timesteps, instabilities lead to a total loss of convergence for ERK4 and ESDC6 and cause EPBM5 to be completely unusable. These limitations make all three exponential integrators only marginally more stable than the fully explicit RK4. Conversely, the IMRK4 method converges across the full range of timesteps and exhibits more than fourth convergence for this problem (approximately 4.5). }
	\label{fig:zds_vanilla_converge_128}	
\end{figure}

To make matters worse, instabilities persists in the ERK4 integrator across a range of coarse timesteps and only disappear when the stepsize has been reduced by a factor of approximately ten. In Figure \ref{fig:zds_vanilla_converge_128} we show convergence diagrams for the ERK4, ESDC6, and EPBM5 integrators. For reference we also include IMRK4 and classical explicit RK4. The reference solution was calculated using RK4 with 200,000 timesteps, and the relative error in all our plots is defined as ${\|\mathbf{y}_{\text{ref}} - \mathbf{y}\|_\infty / \| \mathbf{y}_{\text{ref}} \|}$.

At coarse timesteps the exponential integrators either generate completely inaccurate solutions or are outright unstable. Proper convergence is only achieved when one uses timesteps that are nearly on par with those required to keep explicit RK4 stable. Conversely, the IMRK4 method converges at slightly higher than fourth-order across the full range of timesteps. Lastly, we want to emphasize that this problem presents the same stability challenges for a wide range of exponential integrators. Though we have not completed an exhaustive experiment with all known partitioned exponential integrators, all the families of integrators that we have tested exhibit similar behavior on this problem.

\section{Linear stability}
\label{sec:linear_stability}

In this section we use linear stability analysis \cite[IV.2]{wanner1996solving} to study the stability properties of exponential integrators on non-diffusive problems.  %
 Since we are considering partitioned  integrators, we analyze stability for the non-diffusive, partitioned Dalquist equation
\begin{align}
	y' = i \lambda_1 y + i\lambda_2y, \quad y(0) = y_0, \quad 
	\quad \lambda_1, \lambda_2 \in \mathbb{R}.
	\label{eq:dispersive-dahlquist}
\end{align} 
This equation is a special case of the partitioned Dahlquist equation with $\lambda_1, \lambda_2 \in \mathbb{C}$ that was used to study the general stability of both implicit-explicit \cite{ascher1995implicit,sandu2015generalized,izzo2017highly} and exponential integrators \cite{beylkin1998ELP,cox2002ETDRK4} including ERK4 \cite{krogstad2005IF}, ESDC \cite{buvoli2019esdc} and EPBMs \cite{buvoli2021epbm}. 

We can relate the partitioned Dalquist equation to a general nonlinear system, by observing that (\ref{eq:semilinear_ode}) can be reduced to a decoupled set of partitioned Dahlquist equations provided that $N(t,y)$ is an autonomous, diagonalizable linear operator that shares all its eigenvectors with $\mathbf{L}$. Though this assumption does not hold true for general nonlinear systems, the partitioned Dahlquist equation has nevertheless proven invaluable for studying stability.

When solving (\ref{eq:dispersive-dahlquist}) with a partitioned exponential integrator, we treat the term $i\lambda_1 y$ exponentially and the term $i\lambda_2y$ explicitly. A one-step exponential integrator like ERK or ESDC then reduces to the scalar iteration,
	\begin{align}
		y_{n+1} = R(ik_1,ik_2) y_n	 \quad \text{where} \quad k_1 = h\lambda_1,~ k_2 = h \lambda_2,
		\label{eq:one-step-dahlquist-iteration}
	\end{align}
while a multivalued integrator like EPBM5 reduces to a matrix iteration
\begin{align}
	\mathbf{y}^{[n+1]} &= \mathbf{M}(ik_1,ik_2) \mathbf{y}^{[n]}
\end{align}	
where $\mathbf{M}(ik_1,ik_2) \in \mathbb{C}^{q \times q}$ and $q$ represents the number of input values. The stability region $S$ includes all the $(k_1, k_2)$ pairs that guaranteed a bounded iteration. For the scalar case we simply need a stability function $R$ with magnitude less than one, so that
\begin{align}
	\mathcal{S} &= \left\{ (k_1, k_2) \in \mathbb{R}^2 : |R(ik_1, ik_2)| \le 1 \right\}.
\end{align}	
 For the multivalued case, we require the matrix $\mathbf{M}(ik_1,ik_1)$ to be power bounded. Therefore, we redefine the stability function as $R(ik_1,k_2) = \rho(\mathbf{M}(ik_1,ik_1))$ where $\rho$ denotes the spectral radius. On the boundary, we must take special care to ensure that any eigenvalues of magnitude one are non-defective. 

Since we are only considering real-valued $\lambda_1$ and $\lambda_2$ in (\ref{eq:dispersive-dahlquist}), we can visualize the stability regions for non-diffusive problems using a simple two-dimensional plot. Moreover, all the integrators we consider have stability functions that satisfy $	R(ik_1, ik_2) = R(-ik_1, -ik_2)$.
Therefore, we only need to consider ${k_1 \ge 0}$. This same approach was used in \cite{buvoli2020imex} to visualize the stability of parareal integrators based on implict-explicit Runge-Kutta methods.

\begin{figure}
	\begin{center}
		\begin{tabular}{ccc}
			ERK4 & ESDC6 & EPBM5 \\
			\includegraphics[align=b,width=0.3\linewidth]{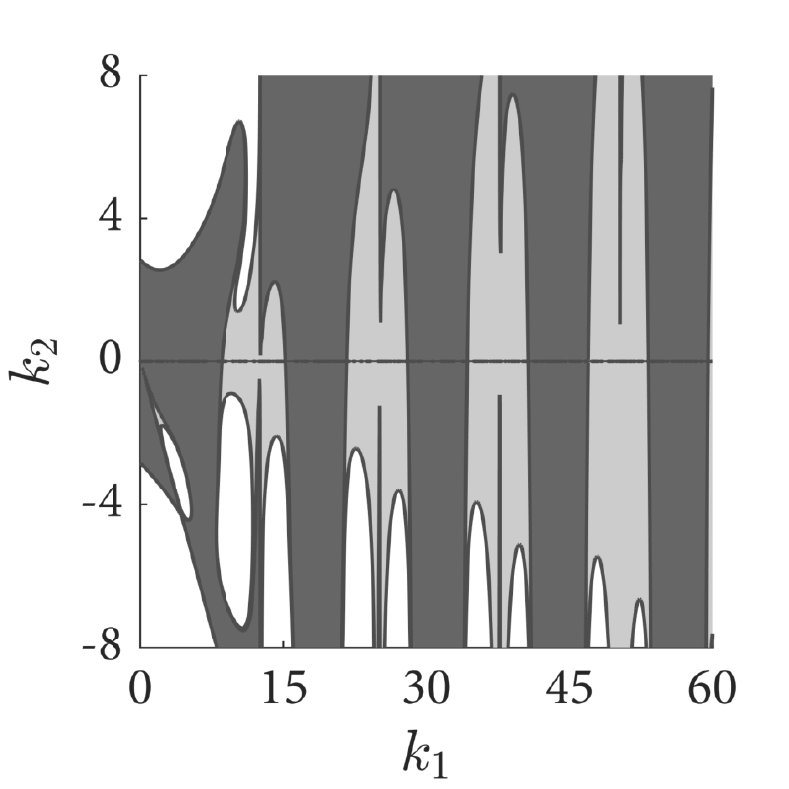} &
			\includegraphics[align=b,width=0.3\linewidth]{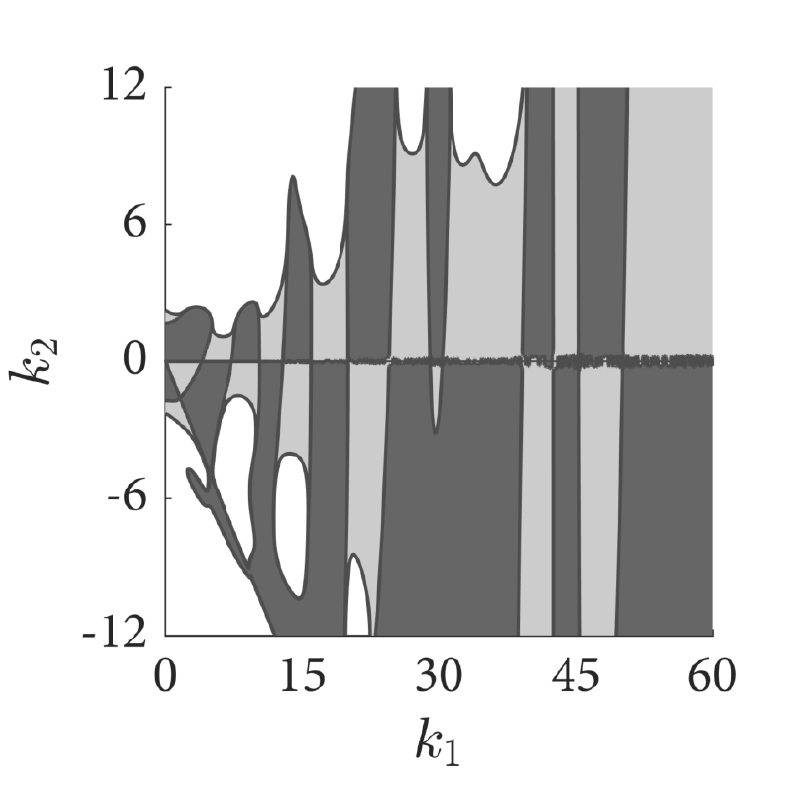} &
			\includegraphics[align=b,width=0.3\linewidth]{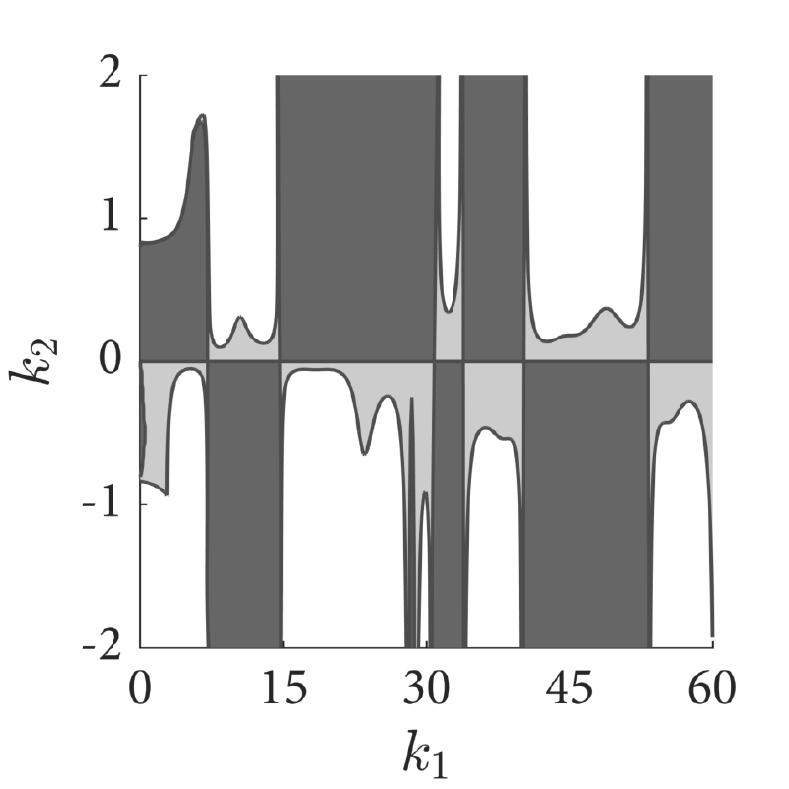}
		\end{tabular}
		\begin{small}
			\textcolor{stability_dark}{\hdashrule[0.0ex]{3em}{5pt}{}} ~~ $|R(ik_1,ik_2)| \le 1$ \hspace{3em}  \textcolor{stability_light}{\hdashrule[0.0ex]{3em}{5pt}{}}  ~~ $|R(ik_1,ik_2)| \le 1.01$	
		\end{small}

		\vspace{1em}
		\begin{tabular}{ccc}
			\includegraphics[align=b,width=0.3\linewidth]{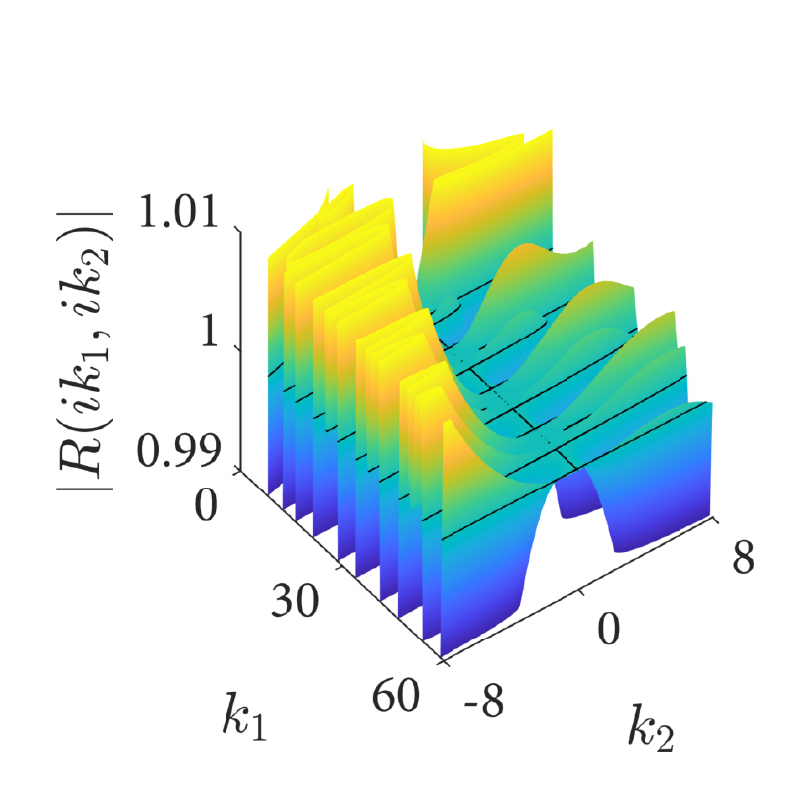} &
			\includegraphics[align=b,width=0.3\linewidth]{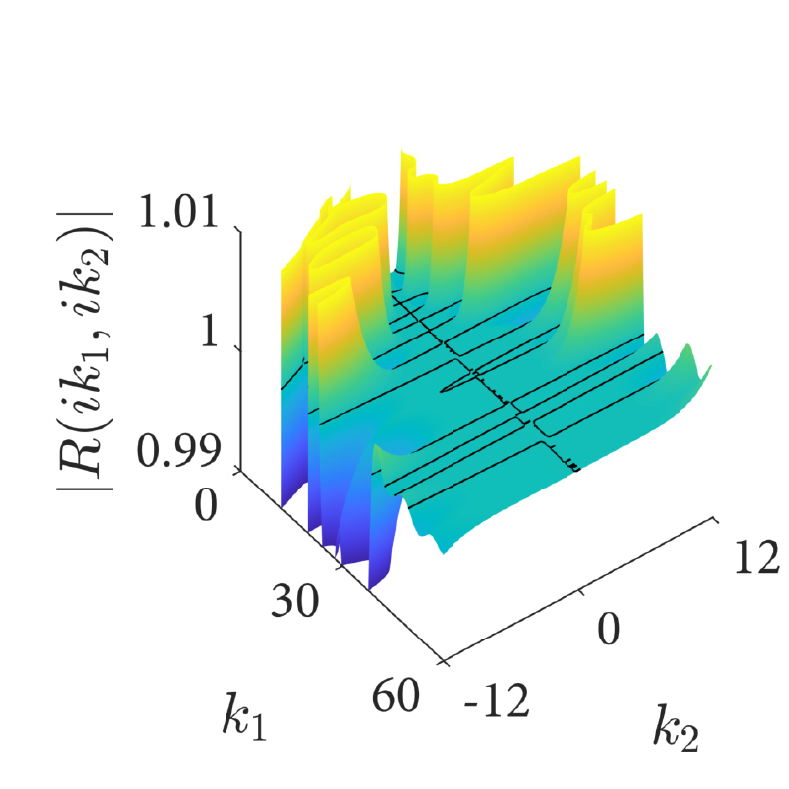} &
			\includegraphics[align=b,width=0.3\linewidth]{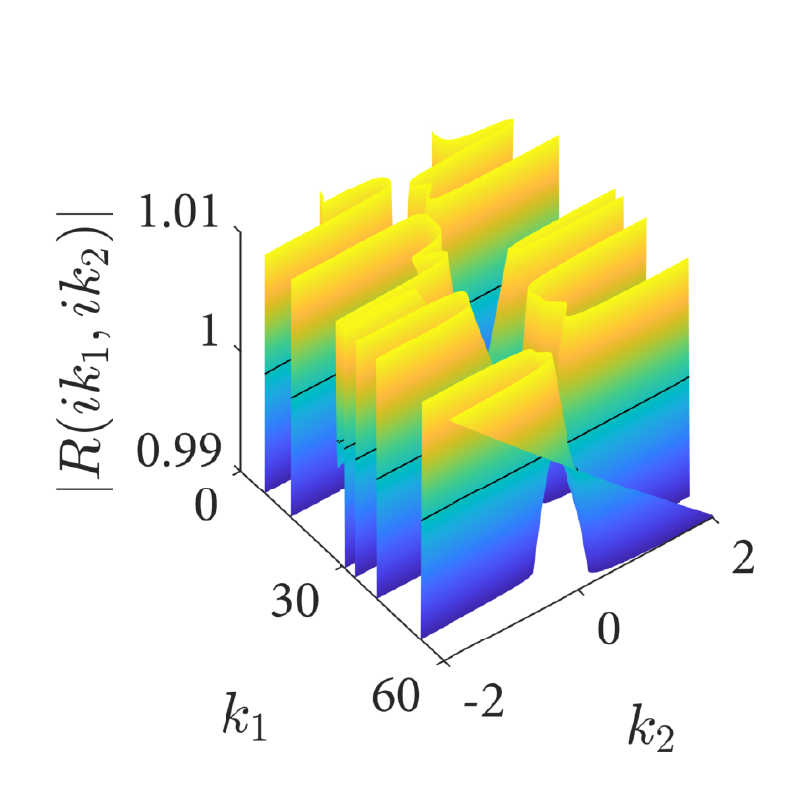}
		\end{tabular}
	\end{center}

	\caption{
		Stability regions $\mathcal{S}$ (top row) and stability functions $R(ik_1,ik_2)$ (bottom row) of the ERK4, ESDC6, and EPBM5 methods for ${0 \le k_1 \le 60}$. To improve readability, the $k_2$ axis is scaled differently for each method. In the top row we show the the proper stability region in dark gray and an extended stability region in light gray that allows for a small amount of instability. 
	}
	\label{fig:vanilla-stability}

\end{figure}

In Figure \ref{fig:vanilla-stability} we show the stability regions $S$ and the stability functions $R(ik_1,ik_2)$ for the ERK4, ESDC6, and EPBM5 methods. If we exclude the line $k_2=0$, which corresponds to $N(t,\mathbf{y})=0$ in (\ref{eq:semilinear_ode}), then the stability regions of all three exponential integrators are disjoint. This characteristic is extremely undesirable and explains the severe timestep restrictions we experienced when solving the ZDS equation. 

Interestingly, the magnitude of the instabilities is very small, and allowing $R(ik_1,ik_2)$ to grow marginally larger than one leads to fully connected regions. However, this observation has little practical significance since we have already seen one example where unmodified exponential integrators produce unusable results. Nevertheless, this phenomenon suggests that exponential integrators will converge successfully on sufficiently small integration intervals that do not allow instabilities to fully manifest. As we will see in Subsection \ref{subsec:kdv-long-time}, this likely explains the existence of successful convergence results found in the literature. %

Another interesting feature of the stability region is the behavior near the line $k_2=0$. Since all exponential integrators treat the linear term exactly, the stability functions all satisfy $R(ik_1,0) = 1$, which pins the stability function right on the border of stability along the $k_2=0$ line. Therefore, even small values of $k_2$ can easily perturb the stability function above one. However, unlike explicit methods, the instabilities decrease in magnitude as $k_1$ increases. This is due to the asymptotic behavior of the exponential integral that approximates the nonlinearity in (\ref{eq:exp_derivation_formula}); namely for any finite $k_2$,
	\begin{align}
		\lim_{k_1 \to \pm \infty}  \left| k_2 \int_{0}^{1} e^{ik_1(s-1)} p(s) ds \right| \sim  \frac{1}{|k_1|} \quad \text{ for any polynomial $p(s)$.}
	\end{align}
	To prove this result one simply needs to apply integration by parts.

\section{Improving stability through repartitioning}
\label{sec:repartitioning}

We now introduce a simple repartitioning strategy to eliminate the stability limitations described in the previous two sections. In short, we propose to repartition the system (\ref{eq:semilinear_ode}) as
\begin{align}
	\mathbf{y}' = \widehat{\mathbf{L}} \mathbf{y} + \widehat{N}(t,\mathbf{y})
	\label{eq:repartioned_semilinear_ode}
\end{align}
where 
\begin{align}
	\widehat{\mathbf{L}} = \mathbf{L} + \epsilon \mathbf{D}, \quad 
	\widehat{N}(t,\mathbf{y}) = N(t,\mathbf{y}) - \epsilon \mathbf{D}\mathbf{y},
	\label{eq:partitioned_operators}
\end{align}
and $\mathbf{D}$ is a diffusive operator. The repartitioned system (\ref{eq:repartioned_semilinear_ode}) is mathematically equivalent to (\ref{eq:semilinear_ode}), however a partitioned exponential integrator now exponentiates $\widehat{\mathbf{L}}$ instead of the original matrix $\mathbf{L}$.

To improve stability, we seek a matrix $\mathbf{D}$ so that the eigenvalues of $\widehat{\mathbf{L}}$ have some small negative real-component. We start by first assuming that $\mathbf{L}$ is diagonalizable so that $\mathbf{L} = \mathbf{U \Lambda U}^{-1}$ and then select
\begin{align}
	\mathbf{D} = -\mathbf{U} \text{abs}(\mathbf{\Lambda}) \mathbf{U}^{-1}.
	\label{eq:D_eigen}
\end{align}
Although this choice may be difficult to apply in practice, it is very convenient for analyzing the stability effects of repartitioning. If we apply (\ref{eq:D_eigen}) to the Dahlquist test problem (\ref{eq:dispersive-dahlquist}), then $\mathbf{D} = -|\lambda_1|$ and we obtain the repartitioned non-diffusive Dahlquist equation
	\begin{align}
		y' = \underbrace{(i \lambda_1  - \epsilon|\lambda_1|)}_{\hat{\mathbf{L}}}y + \underbrace{(i\lambda_2 + \epsilon |\lambda_1|) y}_{\hat{N}(t,y)}.
		\label{eq:repartioned_dispersive_dahlquist}
	\end{align}
The associated stability region for an exponential integrator is 
	\begin{align}
		S &= \left\{ (k_1, k_2) \in \mathbb{R}^2 : \left| \hat{R}(k_1,k_2) \right| \le 1 \right\}, \\
		\hat{R}(k_1,k_2) &= R\left(ik_1 - \epsilon|k_1|, i k_2 + \epsilon |k_1|\right),
		\label{eq:stability_region_damped}	
	\end{align}
where $R$ is the stability function of the unmodified integrator. By choosing
\begin{align}
	\epsilon = 1 / \tan(\tfrac{\pi}{2} + \rho)	\quad \text{for} \quad \rho \in [0, \tfrac{\pi}{2}),
	\label{eq:epsilon_of_rho}
\end{align}
 the single eigenvalue of the partitioned linear operator $\hat{\mathbf{L}}$ %
 is now angled $\rho$ radians off the imaginary axis into the left half-plane. Conversely, the ``nonlinear'' operator $\hat{N}(t,y)$ %
 has been rotated and scaled into the right half-plane. Therefore, in order for the method to stay stable, the exponential functions of $\hat{\mathbf{L}}$ must damp the excitation that was introduced in the nonlinear component.

\begin{figure}[h!]
	
	\begin{center}
	\setlength{\tabcolsep}{0.25em}
	\begin{tabular}{lcccl}
		& ERK4 & ESDC6 & EPBM5 \\
		\rotatebox{90}{ \hspace{2.75em} $\rho = 0$} &
		\includegraphics[align=b,width=0.24\linewidth]{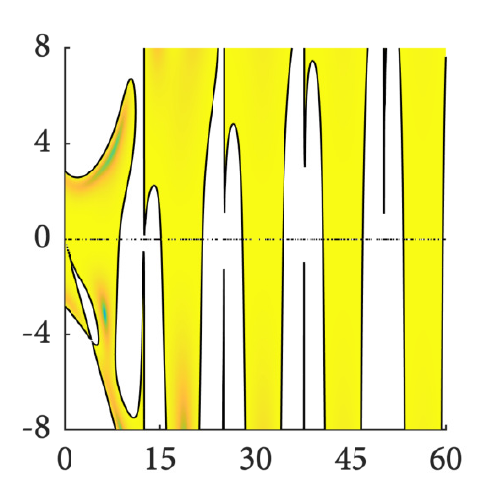} &
		\includegraphics[align=b,width=0.24\linewidth]{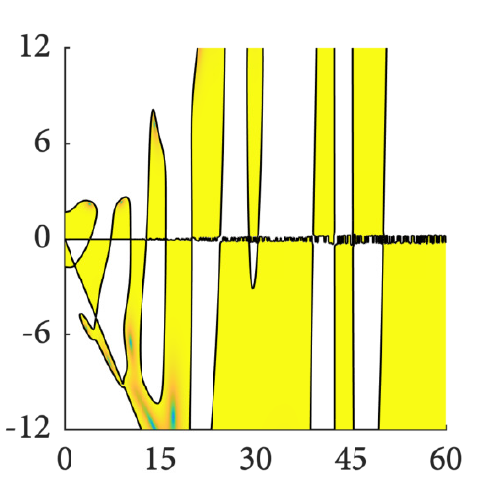} &
		\includegraphics[align=b,width=0.24\linewidth]{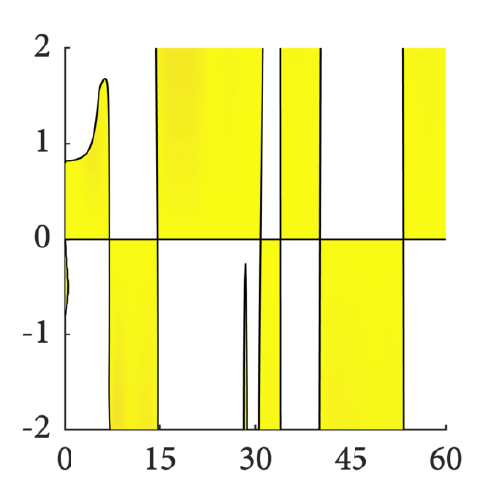}	&
		\includegraphics[align=b,width=0.079\linewidth,trim={100 0 0 0},clip]{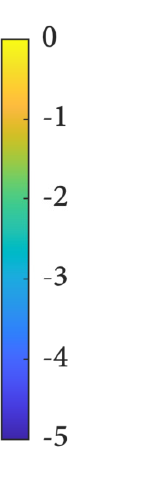} \\
		\rotatebox{90}{ \hspace{2.75em} $\rho = \pi/2048$} &
		\includegraphics[align=b,width=0.24\linewidth]{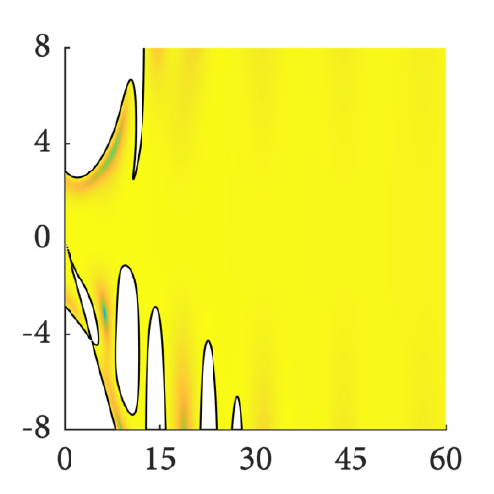} &
		\includegraphics[align=b,width=0.24\linewidth]{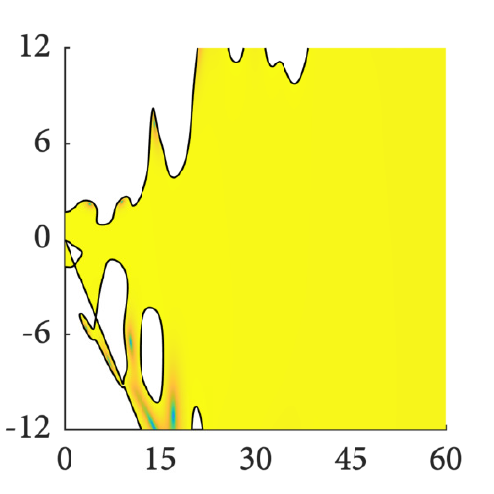} &
		\includegraphics[align=b,width=0.24\linewidth]{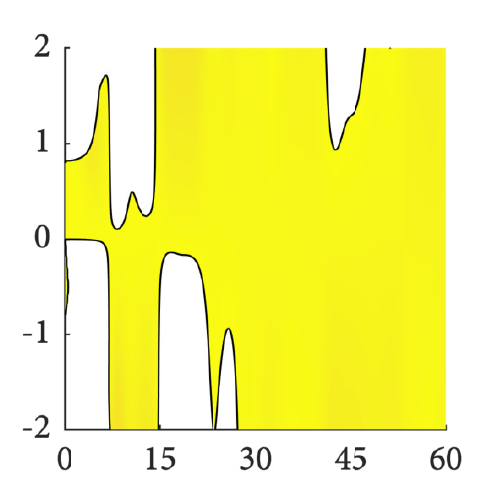} &	 \\
		\rotatebox{90}{ \hspace{2.75em} $\rho = \pi/512$} &
		\includegraphics[align=b,width=0.24\linewidth]{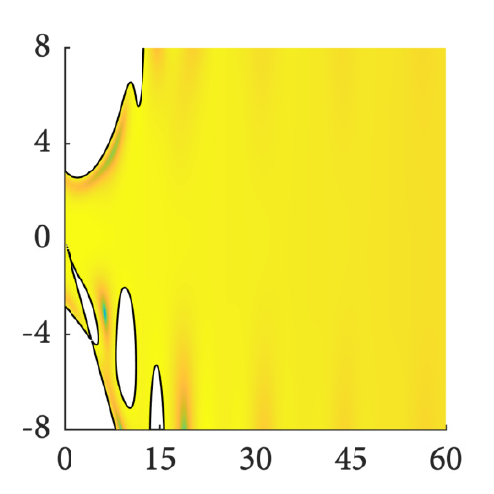} &
		\includegraphics[align=b,width=0.24\linewidth]{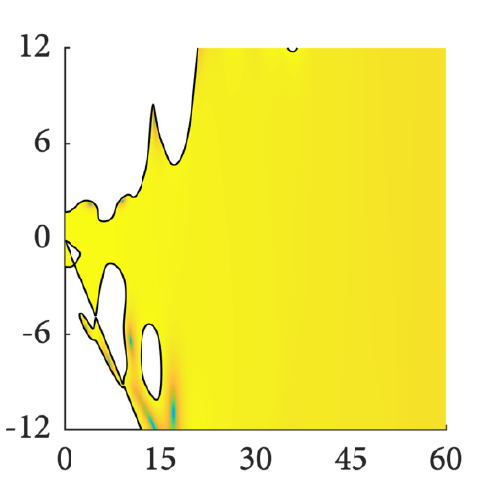} &
		\includegraphics[align=b,width=0.24\linewidth]{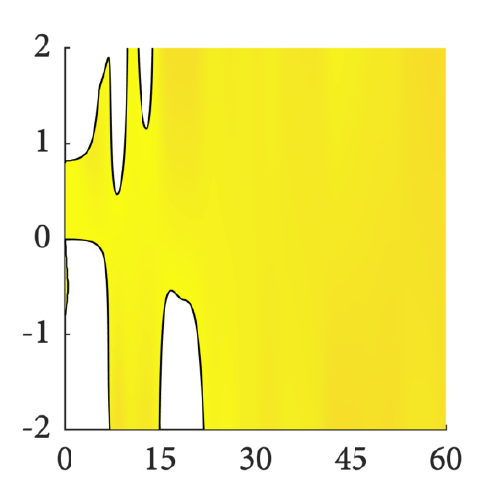} &	\\
		\rotatebox{90}{ \hspace{2.75em} $\rho = \pi/128$} &
		\includegraphics[align=b,width=0.24\linewidth]{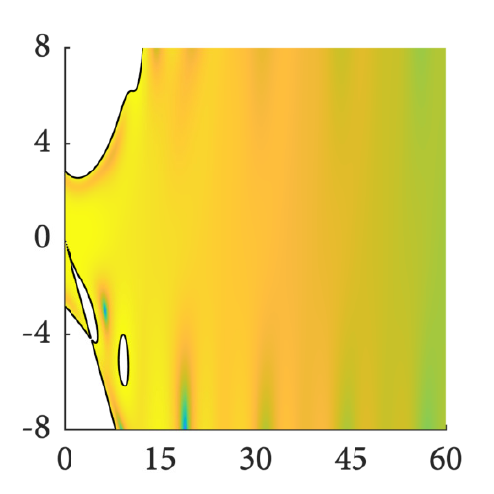} &
		\includegraphics[align=b,width=0.24\linewidth]{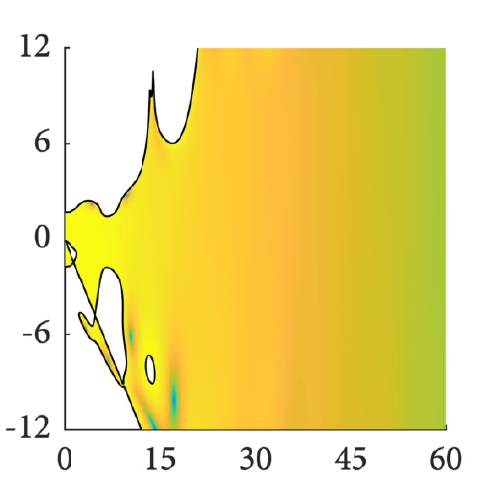} &
		\includegraphics[align=b,width=0.24\linewidth]{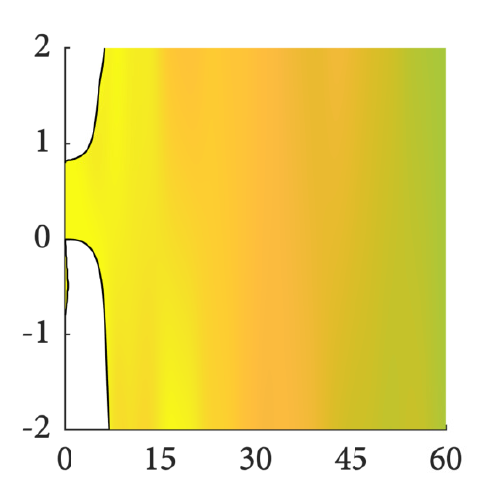} \\
				\rotatebox{90}{ \hspace{2.75em} $\rho = \pi/32$} &
		\includegraphics[align=b,width=0.24\linewidth]{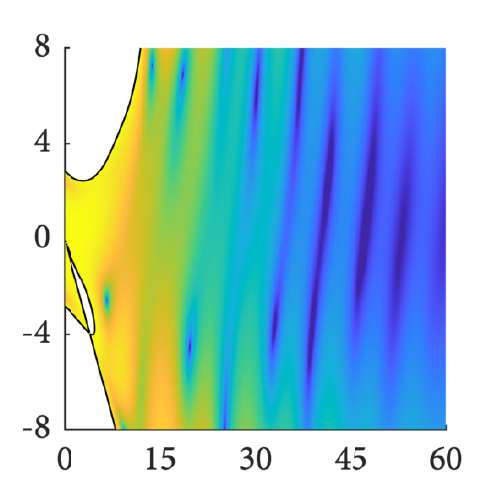} &
		\includegraphics[align=b,width=0.24\linewidth]{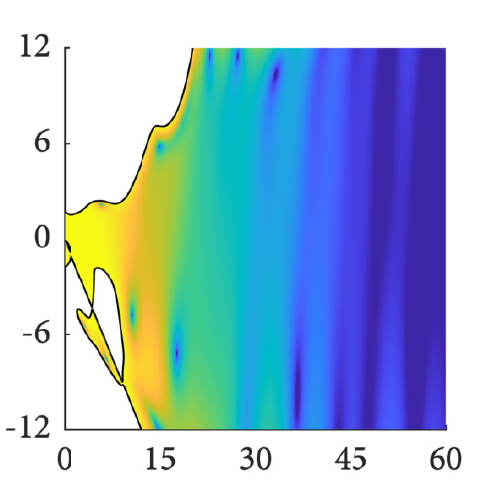} &
		\includegraphics[align=b,width=0.24\linewidth]{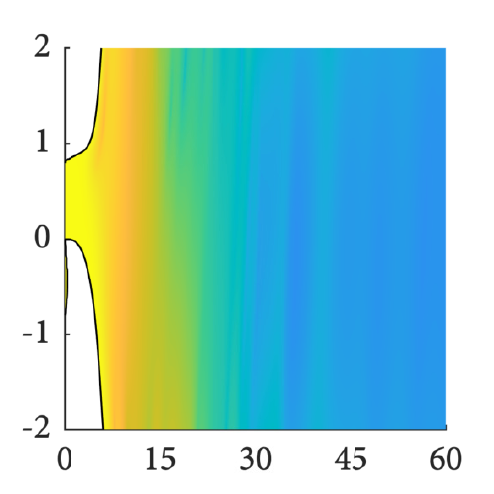} &
	\end{tabular}
	\end{center}
	\vspace{-1em}

	\caption{
		Stability regions $\mathcal{S}$ of the ERK4, ESDC6, and EPBM5 methods for various $\rho$ values. On each plot, the x-axis represents $k_1$ and the y-axis represents $k_2$. The color represents the log of $|\hat{R}(k_1,k_2)|$ from (\ref{eq:stability_region_damped}).
	}
	\label{fig:repartitioned-stability}

\end{figure}

In Figure \ref{fig:repartitioned-stability} we show the stability regions for ERK4, ESDC6, and EPBM5 after repartitioning with different values of $\rho$. Rotating the linear operator by only $\rho = \pi/2048$ radians ($\approx 0.088$ degrees) already leads to large connected linear stability regions for all three methods. This occurs because the magnitude of the partitioned stability function $\hat{R}(k_1,k_2)$ along the line $k_2 = 0$ is now less than one for any $k_1 \ne 0$. In fact, by increasing $\rho$ further, one introduces additional damping for large $k_1$. Therefore, under repartitioning, high-frequency modes will be damped, while low frequency modes will be integrated in a nearly identical fashion to the unmodified integrator with $\rho = 0$. Excluding scenarios where energy conservation is critical, the damping of high-frequency modes is not a serious drawback since large phase errors in the unmodified integrator would still lead to inaccurate solutions (supposing the method is stable in the first place).

\begin{figure}[h!]
	
	\begin{center}	
		\setlength{\tabcolsep}{0.25em}
		\begin{tabular}{ccccr}		
			ERK4 $(\rho = \tfrac{\pi}{4})$ & ESDC6 $(\rho = \tfrac{\pi}{3})$ & EPBM5 $(\rho = \tfrac{\pi}{6})$ \\	
			\includegraphics[width=0.3\linewidth]{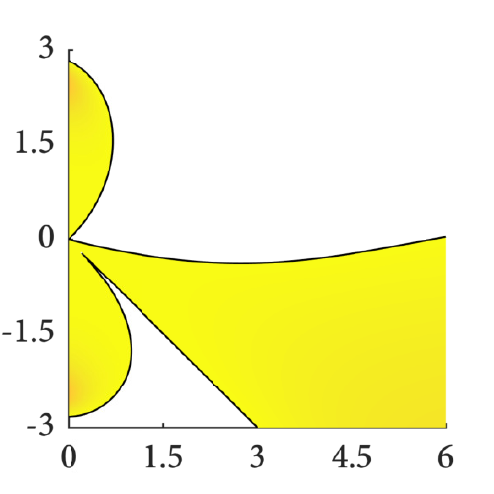} &
			\includegraphics[width=0.3\linewidth]{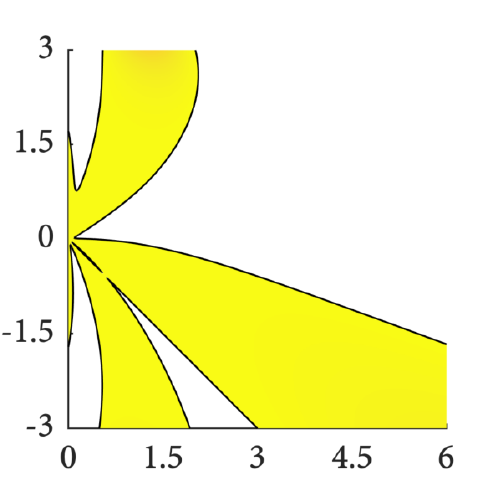} &
			\includegraphics[width=0.3\linewidth]{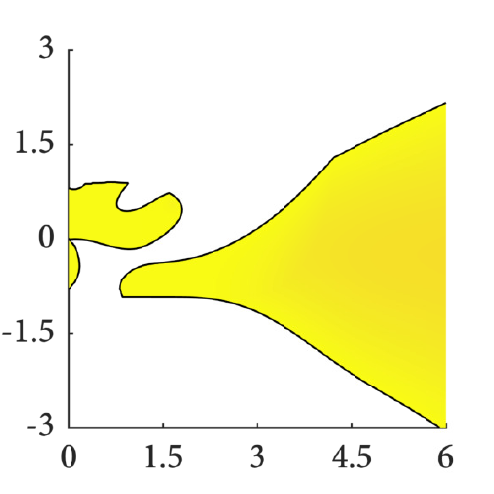} &
			\includegraphics[align=b,width=0.085\linewidth,trim={100 0 0 0},clip]{figures/stability/rho/amp-colorbar}
		\end{tabular}
	\end{center}

	\caption{
		Magnified stability regions that show the effects of adding too much repartitioning. Each integrator family has a different amount of maximal repartitioning before the stability regions split along the $k_2$ line. Amongst the three methods, ESDC6 was the most robust and EPBM5 was the least robust.
	}
	\label{fig:exp-overpartitioning}
	
\end{figure}

Repartitioning cannot be applied indiscriminately and as $\rho$ approaches $\tfrac{\pi}{2}$ one obtains an unstable integrator. To highlight this phenomenon more clearly, we show magnified stability regions in Figure \ref{fig:exp-overpartitioning}, in which we selected sufficiently large $\rho$ values to cause stability region separation for each method. From these results we can see that the maximum amount of allowed repartitioning is integrator dependent, with ESDC6 allowing for the most repartition and EPBM5 the least.

\begin{figure}[h!]
	
	\begin{center}	
		\setlength{\tabcolsep}{0.25em}
			
		\begin{tabular}{ccccr}
		
		$\rho = 0$ & $\rho = \tfrac{\pi}{256}$ & $\rho = \tfrac{\pi}{128}$ & $\rho = \tfrac{\pi}{64}$ \\	
			
		\includegraphics[width=0.23\linewidth]{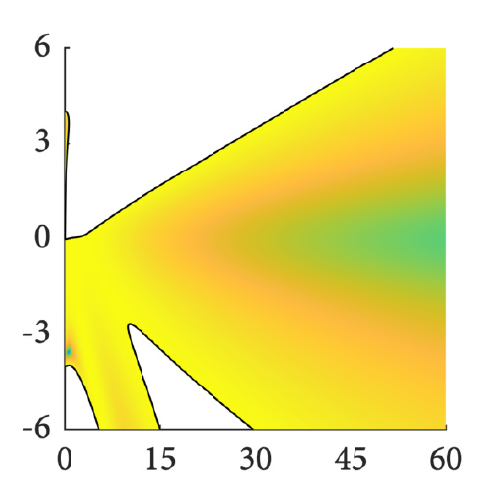} &
		\includegraphics[width=0.23\linewidth]{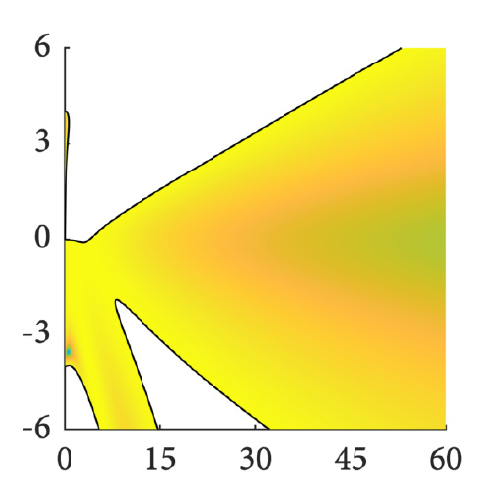} &
		\includegraphics[width=0.23\linewidth]{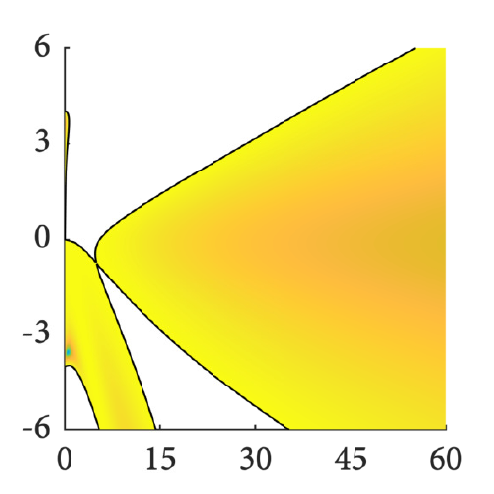} &
		\includegraphics[width=0.23\linewidth]{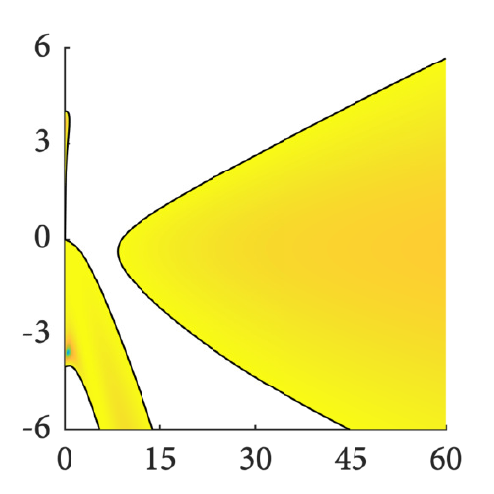} &
		\includegraphics[align=b,width=0.063\linewidth,trim={100 0 0 0},clip]{figures/stability/rho/amp-colorbar}
			
		\end{tabular}
	\end{center}

	\caption{
		Stability regions for the repartitioned IMRK4 methods with four values of $\rho$. For even small values of $\rho$, instabilities form along the $k_2 = 0$ line near the origin, and for large values of $\rho$, the stability regions separates. The effect appears for even smaller values of $\rho$, however we did not include these plots here, since the effect is not visible at  this level of magnification.	
	}
	\label{fig:imex-partitioning}
	
\end{figure}

Finally, we note that this repartitioning technique can also be applied to implicit-explicit methods. However, on all the methods we tried, we found that the repartitioning rapidly destabilizes the integrator. In Figure \ref{fig:imex-partitioning} we present the stability regions for IMRK4 using different $\rho$ values and show that stability along the $k_2=0$ line is lost even for small $\rho$ values. The stability region corresponding to $\rho = 0$ can be compared with the stability regions of the exponential integrators in Figure \ref{fig:repartitioned-stability} to see how the damping properties of a repartitioned exponential methods compare to those of an IMEX method.

\clearpage
\subsection{Solving ZDS with repartitioning}

We now validate our stability results by solving the ZDS equation (\ref{eq:zds}) using several different choices for the diffusive operator $\mathbf{D}$. Since we are solving in Fourier space where $\mathbf{L} = \text{diag}(i \mathbf{k}^3)$, we can easily implement (\ref{eq:D_eigen}) by selecting $\mathbf{D} = -\text{diag}(|\mathbf{k}|^3)$. However, for more general problems, it will often not be possible to compute the eigenvalue decomposition of the operator $\mathbf{L}$. Therefore, to develop a more practical repartitioning, we consider a generic, repartitioned, semilinear partial differential equation
	\begin{align}
		u_t = \underbrace{L[u] + \epsilon D[u]}_{\hat{L}[u]} + \underbrace{N(t,u) - \epsilon D[u]}_{\hat{N}(t,u)}
		\label{eq:repartitioned-semilinear-pde}
	\end{align}
	where the spatial discretizations of $\hat{L}[u]$ and $\hat{N}(t,u)$ become the $\hat{\mathbf{L}}$ and $\hat{N}(t, \mathbf{y})$ in (\ref{eq:partitioned_operators}), and the linear operator $D[u]$ is the continuous equivalent of the matrix $\mathbf{D}$. A natural choice for a diffusive operator in one spatial dimension is the even-powered derivative
		\begin{align}
			D[u] = 
			\begin{cases}
 				\frac{\partial^k u}{\partial x^k} & k \equiv 0\pmod{4}, \\
 				-\frac{\partial^k u}{\partial x^k} & k \equiv 2\pmod{4}. \\
 			\end{cases}
		\end{align}
	This operator can be easily implemented for different spatial discretizations and boundary conditions. Moreover, it can be generalized to higher dimensional PDEs by adding partial derivatives in other dimensions. 
	
\begin{figure}[h!]
	
	\centering
	
	{\small \hspace{4em} $\mathbf{D} = -\text{diag}\left(|\mathbf{k}|^3\right)$}
	
	\includegraphics[width=\linewidth,trim={20 0 52 0},clip]{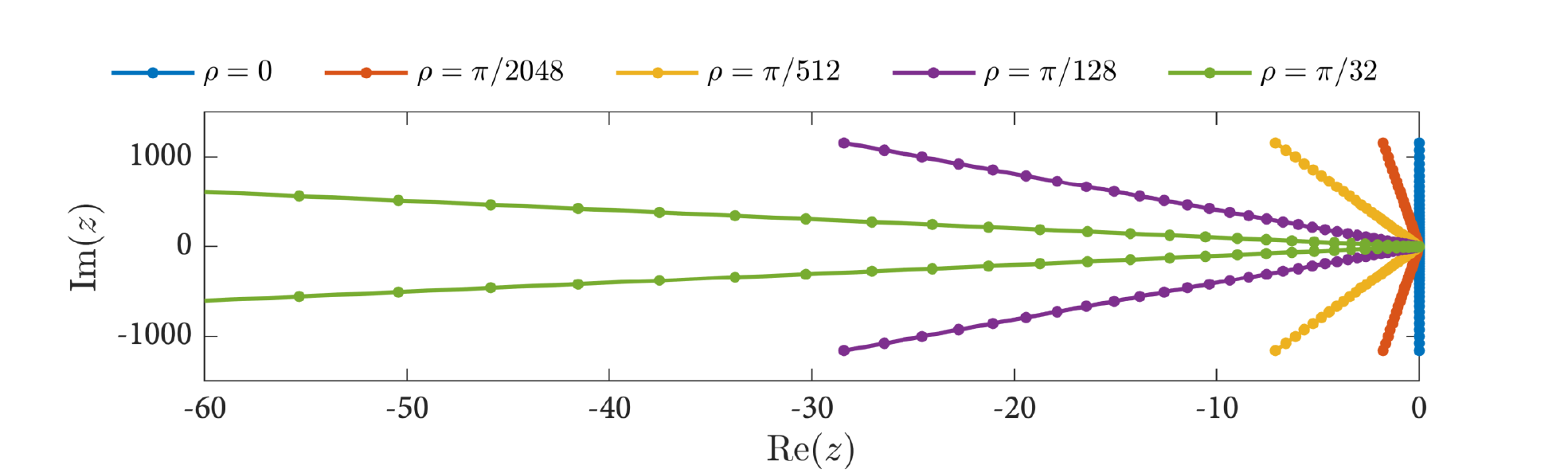}
	
	\vspace{1em}
	{\small \hspace{4em} $\mathbf{D} = -\text{diag}\left(\mathbf{k}^2\right)$}
	
	\includegraphics[width=\linewidth,trim={20 0 52 0},clip]{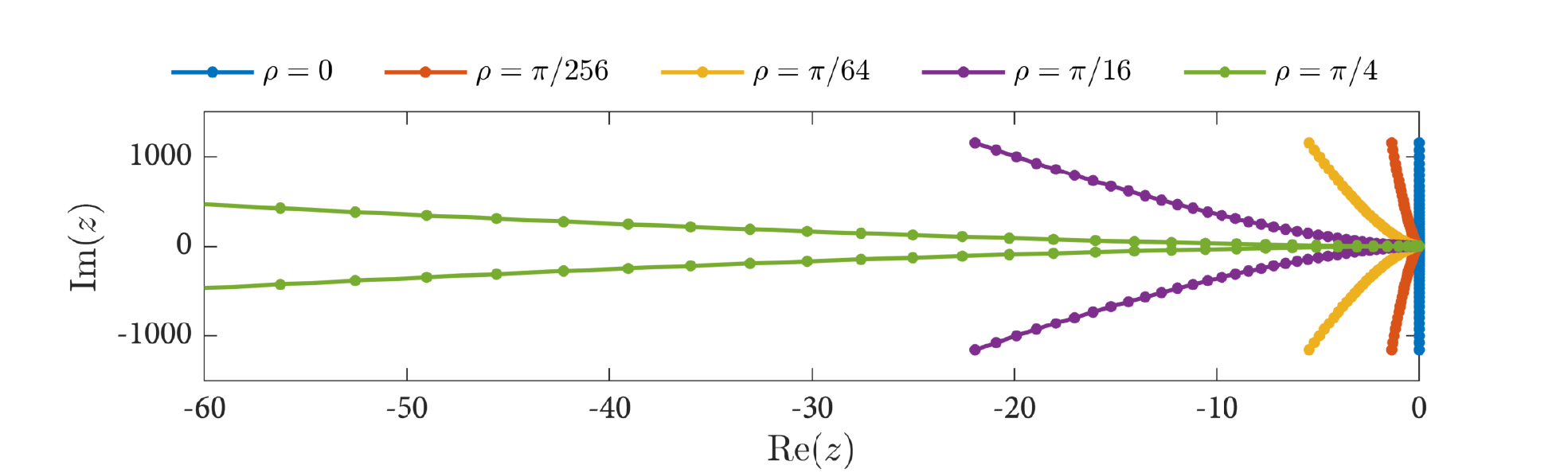}
	
	\vspace{1em}
	{\small \hspace{4em} $\mathbf{D} = -\mathbf{I}$}
	
	\includegraphics[width=\linewidth,trim={20 0 52 0},clip]{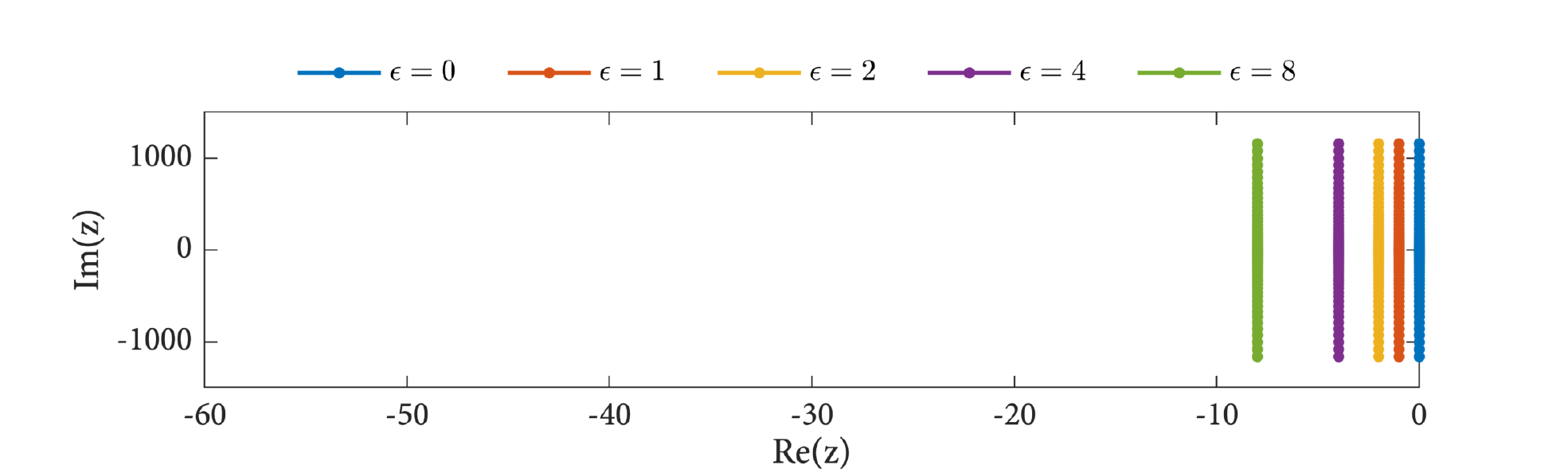}

	\caption{Spectrum of the repartitioned linear operator $\hat{\mathbf{L}} = \mathbf{L} + \epsilon \mathbf{D}$, for three choices of $\mathbf{D}$ and multiple $\epsilon$. In the first two plots $\epsilon$ is selected according to (\ref{eq:epsilon_of_rho}). The choices of $\rho$ for the second-order repartitioning where selected to achieve comparable damping to third-order repartitioning.}%

	\label{fig:zds_partitioned_linear_operator}
\end{figure}

	The only remaining question is how to choose $k$. To avoid increasing the number of required boundary conditions, it is advantageous if $k$ is smaller than or equal to the highest derivative order found in $L[u]$.  Therefore, in addition to (\ref{eq:D_eigen}), we also consider two additional repartitionings for the ZDS equation that are based on the zeroth and the second spatial derivative of $u(x,t)$. Below, we describe each repartitioning in detail, and in Figure \ref{fig:zds_partitioned_linear_operator} we also plot the spectrum of the corresponding repartitioned linear operators $\hat{\mathbf{L}}$.
	\begin{itemize}
		\item {\bf Third-order repartitioning}. The diffusive operators are
				\begin{align}
					D[u] &= \mathcal{F}^{-1}(|k|^3) \ast u &&\text{(Continuous  -- physical space),}& \\
					\mathbf{D} &= -\text{diag}(|\mathbf{k}|^3)	 &&\text{(Discrete -- Fourier space),} &
				\end{align}		
		where $\mathcal{F}^{-1}$ denotes the inverse Fourier transform and $\ast$ is a convolution. This choice is equivalent to (\ref{eq:D_eigen}). We choose $\epsilon$ according to (\ref{eq:stability_region_damped}), so that the eigenvalues of $\hat{\mathbf{L}} = \mathbf{L} + \epsilon \mathbf{D}$ lie on the curve 
			\begin{align*}
				r e^{i(\pi/2 + \rho)} \cup r e^{i(3 \pi/2 - \rho)} \quad (r \ge 0).
			\end{align*}
		For this repartitioning we select $\epsilon$ using $\rho = \frac{\pi}{2048}$, $\frac{\pi}{512}$, $\frac{\pi}{128}$, and $\frac{\pi}{32}$.

		\item {\bf Second-order repartitioning}. The diffusive operators are
			\begin{align}
				D[u] &= u_{xx} &&\text{(Continuous  -- physical space),}\\
				\mathbf{D} &= -\text{diag}(\mathbf{k}^2) &&\text{(Discrete  -- Fourier space),} & 	
			\end{align} 	
			and we again choose $\epsilon$ according to (\ref{eq:stability_region_damped}). Compared to the previous choice, second-order repartitioning over-rotates eigenvalues with magnitude less then one, and under-rotates eigenvalues with magnitude larger than one. Therefore we require larger $\rho$ values to achieve similar damping effects as third-order repartitioning (see Figure \ref{fig:zds_partitioned_linear_operator}); in particular we select $\epsilon$ using $\rho = \frac{\pi}{256}$, $\frac{\pi}{64}$, $\frac{\pi}{16}$, and $\frac{\pi}{4}$.
		
		\item {\bf Zeroth-order repartitioning}. The diffusive operators are
			\begin{align}
				D[u] &= -u &&\text{(Continuous -- physical space),}\\
				\mathbf{D} &= -\mathbf{I} &&\text{(Discrete  -- Fourier space).}	
			\end{align} 
		This choice translates every eigenvalue of the linear operator $\mathbf{L}$ by a fixed amount $\epsilon$ into the left-hand plane. We consider ${\epsilon = 1, 2, 4}$ and $8$.
		
	\end{itemize}
	
\begin{figure}[h!]
	
	\begin{center}
		\setlength{\tabcolsep}{0.1em}
		\begin{tabular}{llcccc}
			& \hspace{.5em} & $\mathbf{D} = -\text{diag}\left(|\mathbf{k}|^3\right)$ & $\mathbf{D} = -\text{diag}\left(\mathbf{k}^2\right)$ & $\mathbf{D} = -\mathbf{I}$ \\
			\rotatebox{90}{ \hspace{5.9em} ERK4 } & &
			\includegraphics[align=b,width=0.30\linewidth]{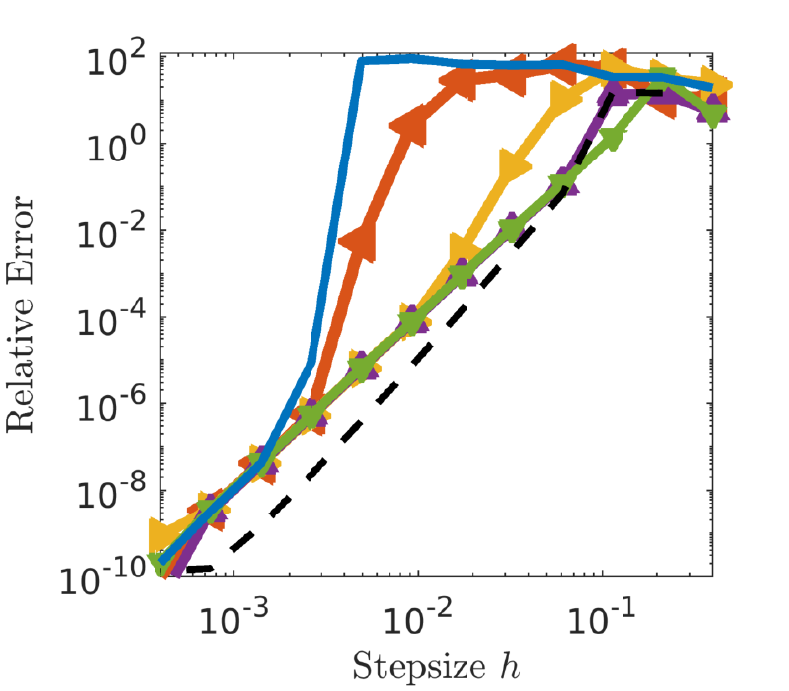} & 
			\includegraphics[align=b,width=0.30\linewidth]{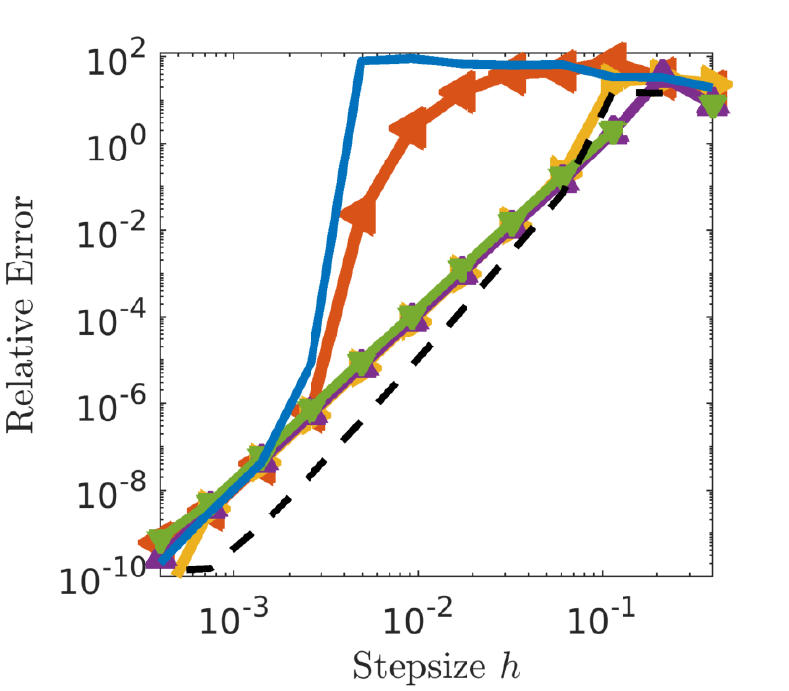} & 
			\includegraphics[align=b,width=0.30\linewidth]{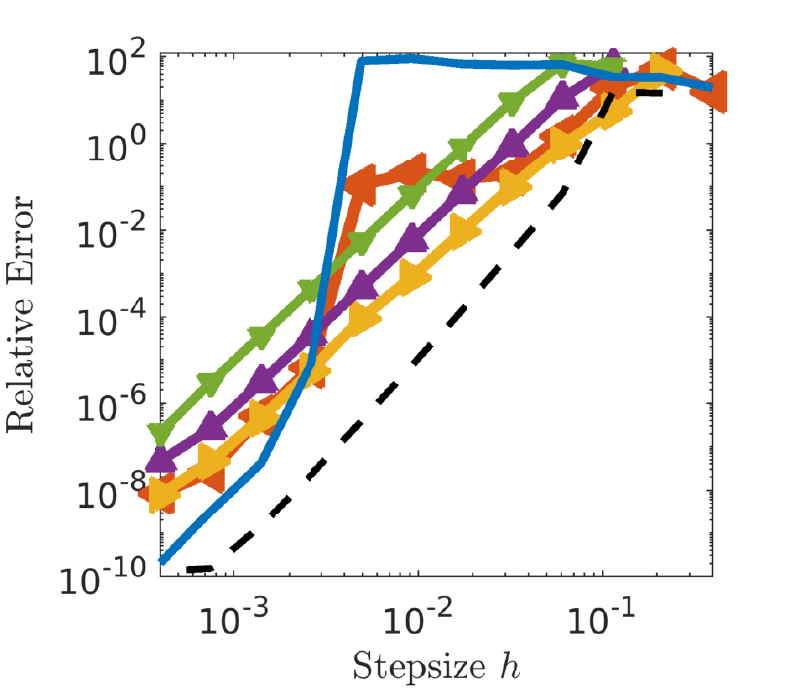} \\
			\rotatebox{90}{ \hspace{5.6em} ESDC6 } & &
			\includegraphics[align=b,width=0.30\linewidth]{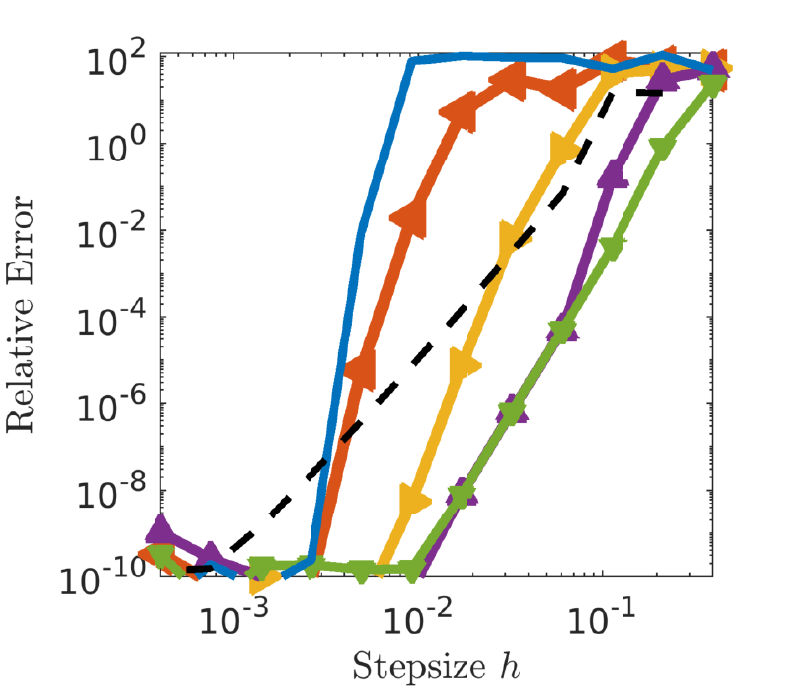} &
			\includegraphics[align=b,width=0.30\linewidth]{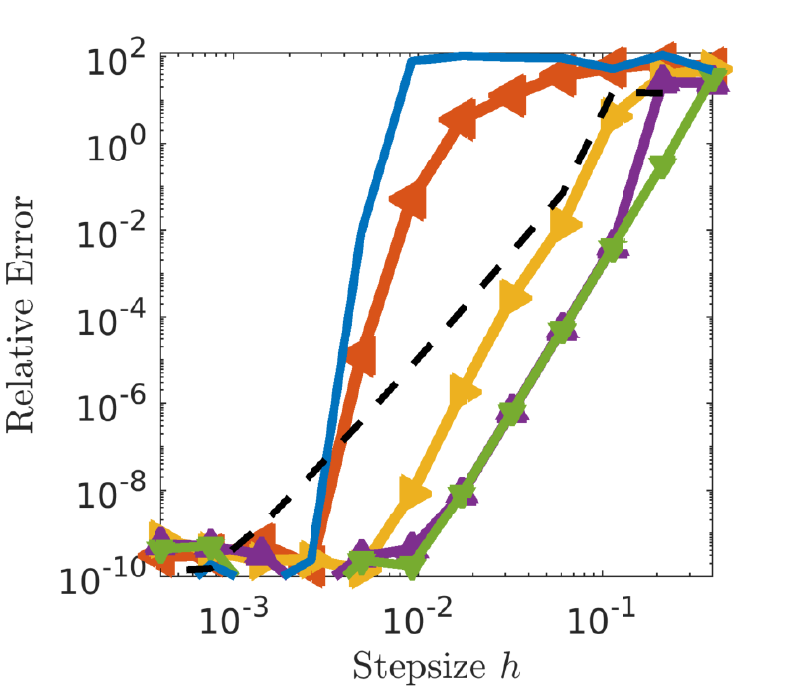} & 
			\includegraphics[align=b,width=0.30\linewidth]{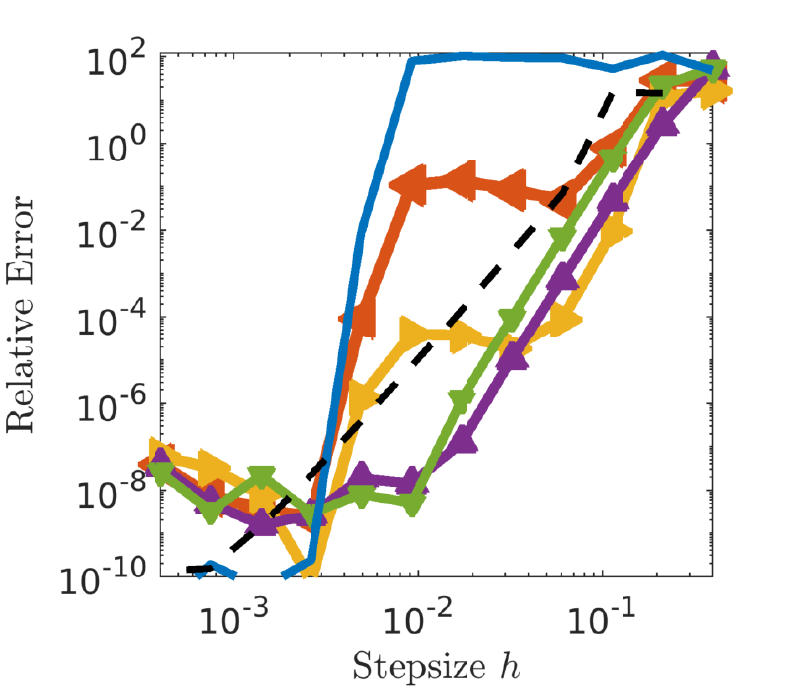} \\
			\rotatebox{90}{ \hspace{5.45em} EPBM5 } & &
			\includegraphics[align=b,width=0.30\linewidth]{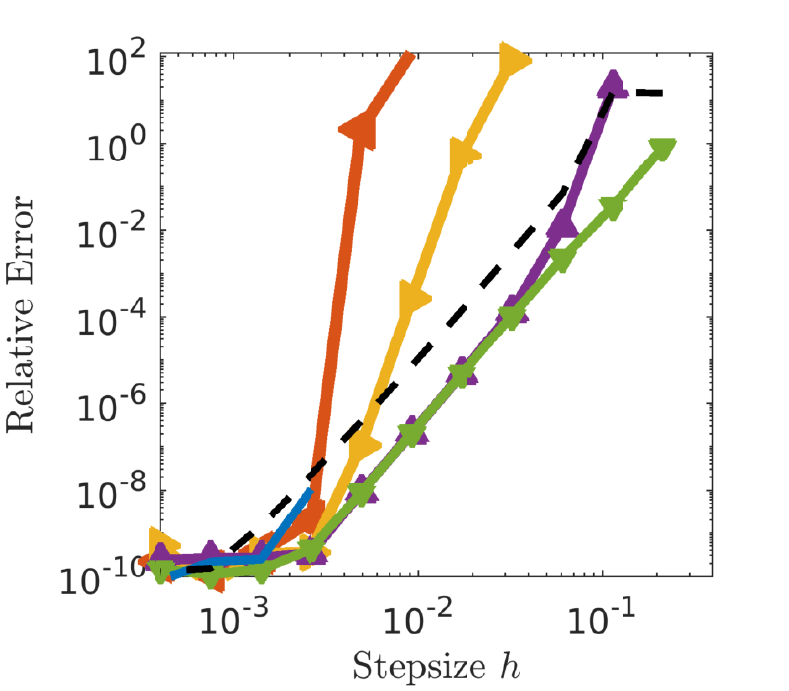} &
			\includegraphics[align=b,width=0.30\linewidth]{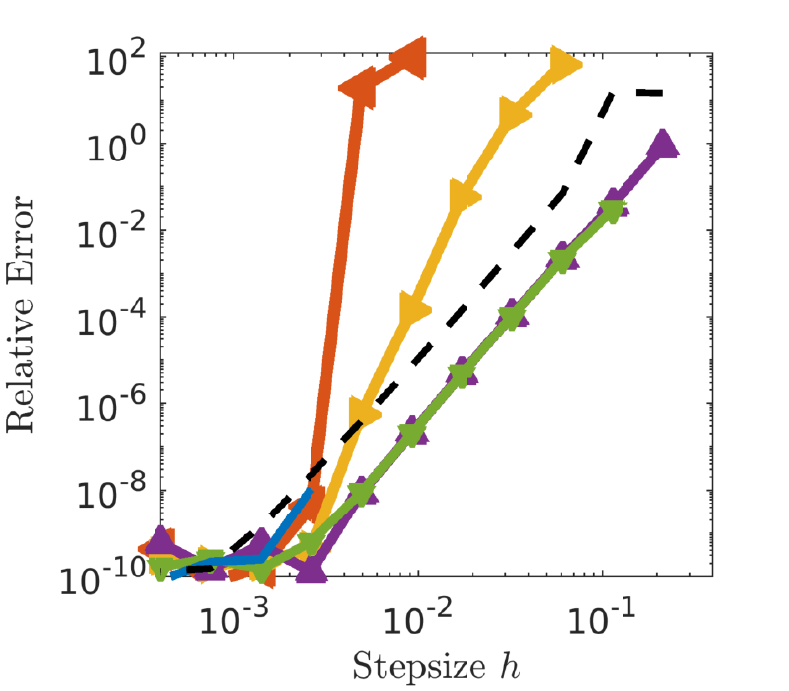} &
			\includegraphics[align=b,width=0.30\linewidth]{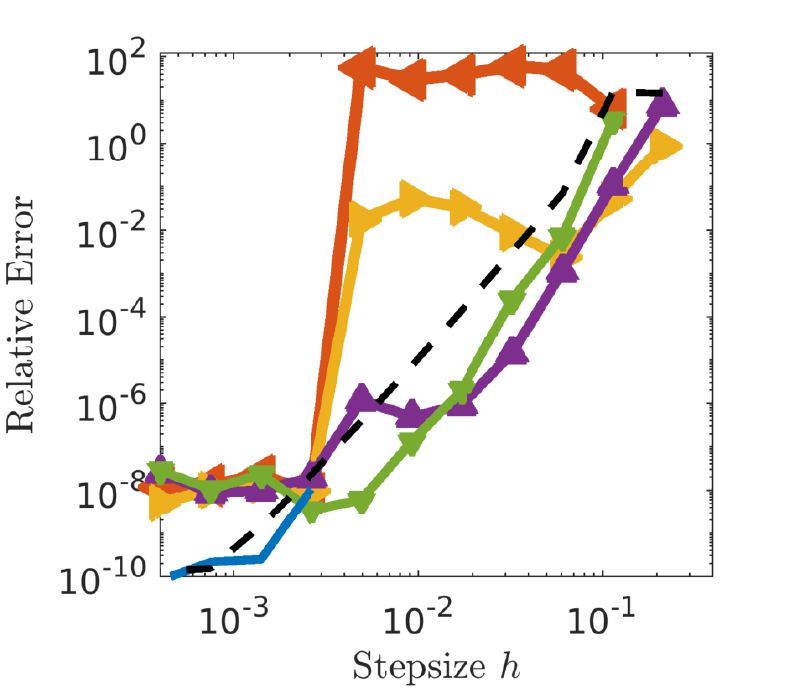}
		\end{tabular}
		
		\begin{small}
		\hspace{3.5em}
		\begin{tabular}{c|c|c}
		\begin{tabular}{ll}
			\includegraphics[align=m,width=0.7cm,trim={560 190 115 14},clip]{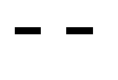} IMRK4 ~~~~ & 
			\includegraphics[align=m,width=0.7cm,trim={100 190 580 14},clip]{figures/experiment/legend-rotation} $\rho = 0$ \\			
			\includegraphics[align=m,width=0.7cm,trim={180 190 505 14},clip]{figures/experiment/legend-rotation} $\rho = \frac{\pi}{2048}$ &	
			\includegraphics[align=m,width=0.7cm,trim={280 190 400 14},clip]{figures/experiment/legend-rotation} $\rho = \frac{\pi}{512}$ ~~ \\ 
			\includegraphics[align=m,width=0.7cm,trim={375 190 305 14},clip]{figures/experiment/legend-rotation} $\rho = \frac{\pi}{128}$ ~~ &
			\includegraphics[align=m,width=0.7cm,trim={475 190 205 14},clip]{figures/experiment/legend-rotation} $\rho = \frac{\pi}{32}$
		\end{tabular} \hspace{0.75em}
		&
		\hspace{0.75em}
		\begin{tabular}{ll}
			\includegraphics[align=m,width=0.7cm,trim={560 190 115 14},clip]{figures/experiment/legend-rotation} IMRK4 ~~~~ & 
			\includegraphics[align=m,width=0.7cm,trim={100 190 580 14},clip]{figures/experiment/legend-rotation} $\rho = 0$ \\			
			\includegraphics[align=m,width=0.7cm,trim={180 190 505 14},clip]{figures/experiment/legend-rotation} $\rho = \frac{\pi}{512}$ &	
			\includegraphics[align=m,width=0.7cm,trim={280 190 400 14},clip]{figures/experiment/legend-rotation} $\rho = \frac{\pi}{128}$ ~~ \\ 
			\includegraphics[align=m,width=0.7cm,trim={375 190 305 14},clip]{figures/experiment/legend-rotation} $\rho = \frac{\pi}{32}$ ~~ &
			\includegraphics[align=m,width=0.7cm,trim={475 190 205 14},clip]{figures/experiment/legend-rotation} $\rho = \frac{\pi}{4}$
		\end{tabular} \hspace{0.75em}
		&
		\hspace{0.75em}
		\begin{tabular}{ll}
			\includegraphics[align=m,width=0.7cm,trim={560 190 115 14},clip]{figures/experiment/legend-rotation} IMRK4 ~~~~ & 
			\includegraphics[align=m,width=0.7cm,trim={100 190 580 14},clip]{figures/experiment/legend-rotation} $\epsilon = 0$ \\			
			\includegraphics[align=m,width=0.7cm,trim={180 190 505 14},clip]{figures/experiment/legend-rotation} $\epsilon = 1$ &	
			\includegraphics[align=m,width=0.7cm,trim={280 190 400 14},clip]{figures/experiment/legend-rotation} $\epsilon = 2$ ~~ \\ 
			\includegraphics[align=m,width=0.7cm,trim={375 190 305 14},clip]{figures/experiment/legend-rotation} $\epsilon = 4$ ~~ &
			\includegraphics[align=m,width=0.7cm,trim={475 190 205 14},clip]{figures/experiment/legend-rotation} $\epsilon = 8$
		\end{tabular} \hspace{1.2em}
		\end{tabular}
		\end{small}		
	\end{center}
	\vspace{-1em}
	
	\caption{
		Convergence diagrams for the repartitioned ERK4, ESDC6, and EPBM5 integrators on the ZDS equation using third, second, and zeroth order repartitioning. For reference we also include the unpartitioned IMRK4 integrator in each plot. The columns correspond to different choices for the matrix $\mathbf{D}$, and each row to different integrators. Color denotes the value of $\rho$ from (\ref{eq:epsilon_of_rho}) for the first two columns and the value of $\epsilon$ for the third column.  
	}
	\label{fig:repartitioned-stability-a}
	
\end{figure}

In Figures \ref{fig:repartitioned-stability-a} we present convergence diagrams for ERK4, ESDC6, and EPBM5 using each of the three repartitioning strategies. Overall, repartitioning resolved the stability issues and enabled the use of exponential integrators for efficiently solving the ZDS equation. We summarize the results as follows.

\vspace{1em}
{\em \noindent Third-order repartitioning.}
Adding even a small amount of third-order repartitioning immediately improves the convergence properties of the exponential integrators. For $\rho = \pi/128$, all integrators achieve proper convergence across the full range of stable timesteps. Moreover, adding additional repartitioning does not damage the accuracy so long as the underlying method remains stable.

\vspace{1em}
{\em \noindent Second-order repartitioning.}
Second-order repartitioning is able to achieve nearly identical results to third-order repartitioning, provided that larger $\rho$-values are used. %
Overall, the results are not surprising since the spectrums of the corresponding linear operators shown in Figure \ref{fig:zds_partitioned_linear_operator} look very similar. The main disadvantage of second-order repartitioning is that $\rho$ needs to tuned to ensure that the highest modes have been sufficiently rotated.
	
\vspace{1em}
{\em \noindent Zeroth-order repartitioning.} Zeroth order repartitioning is extremely simple to implement, however it is also the least effective at improving convergence and preserving accuracy. A small $\epsilon$ does not introduce enough damping and the convergence curves are improved but not fully restored. On the other hand, large $\epsilon$ values stabilize stiff modes, however since all the eigenvalues are shifted by an equal amount, the repartitioning damages the accuracy of non-stiff modes. This leads to convergence curves that have been shifted to the left since we have effectively worsened the error constant of the exponential integrator. Zeroth-order repartitioning also negatively impacted the sensitivity of the integrator to roundoff errors, and we were unable to obtain the solution with a relative error of less than approximately $10^{-8}$.
	
\clearpage

\subsection{Comparing repartitioning to hyperviscosity}

Time dependent equations with no diffusion are known to cause stability issues for numerical methods, and a commonly applied strategy is to add a hyperviscous term to the equation right-hand-side (see for example \cite{jablonowski2011pros,Ullrich2018-fs}). To avoid destroying the convergence properties of an integrator of order $q$, the magnitude of this term is typically proportional to the stepsize of the integrator $\Delta t$ raised to the power of $q+1$.

In the context of the continuous semilinear partial differential equation
\begin{align}
	u_t = L[u] + N(t,u)
	\label{eq:generic-semilinear-pde}	
\end{align}
this is equivalent to considering a new equation with a vanishing diffusive operator $\tilde{D}[u]$ added to the right-hand-side so that
	\begin{align}
		u_t = L[u] + (\Delta  t)^{q+1} \gamma \tilde{D}[u] + N(t,u),
		\label{eq:hyperviscous-semilinear-pde}
	\end{align}
	where $\gamma$ is a constant that controls the strength of the diffusion. One then approximates the solution to (\ref{eq:generic-semilinear-pde}) by numerically integrating (\ref{eq:hyperviscous-semilinear-pde}). The improvement to stability comes from the fact that we have replaced the original discretized linear operator $\mathbf{L}$ with ${\tilde{\mathbf{L}} = \mathbf{L} + (\Delta t)^{q+1}\gamma \tilde{\mathbf{D}}}$. 
	
	Unlike repartitioning, we are no longer adding and subtracting the new operator. We must therefore ensure that $D[u]$ does not damage the accuracy of slow modes as they typically contain the majority of useful information. For this reason, $D[u]$ is generally chosen to be a high-order even derivative since these operators have a negligible effect on low frequencies while causing significant damping of high frequencies. %
	
\begin{figure}[h]
	
	\begin{center}
		
		\begin{footnotesize}
		\setlength{\tabcolsep}{0.1em}
		\begin{tabular}{ccc}
			$\tilde{\mathbf{D}} = -\text{diag}(\mathbf{k}^4)$, $\gamma=10^4 \omega$ & $\tilde{\mathbf{D}} = - \text{diag}(\mathbf{k}^6)$, $\gamma=10^2\omega$ & $\tilde{\mathbf{D}} = -\text{diag}(\mathbf{k}^8)$, $\gamma = \omega$ \\
			\includegraphics[align=b,width=0.33\linewidth]{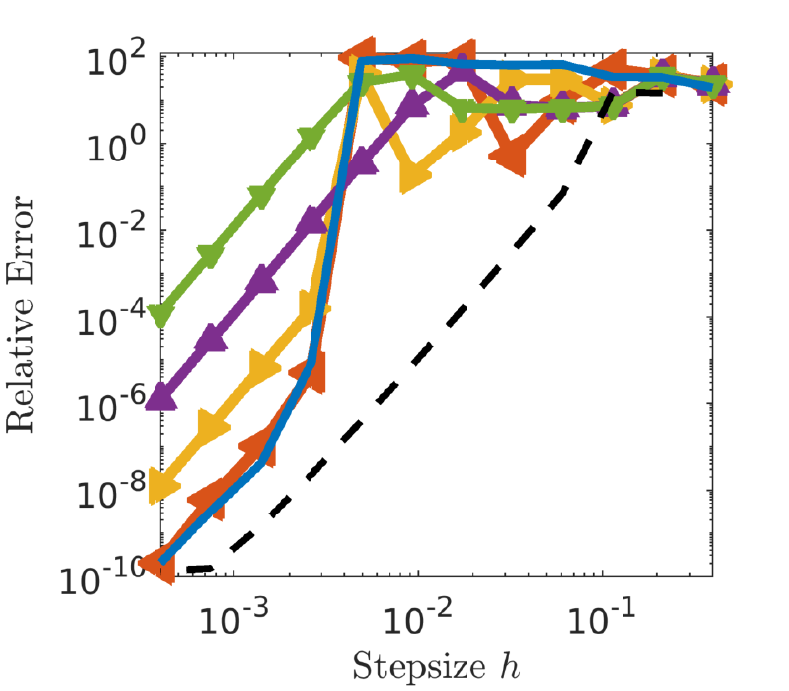} &
			\includegraphics[align=b,width=0.33\linewidth]{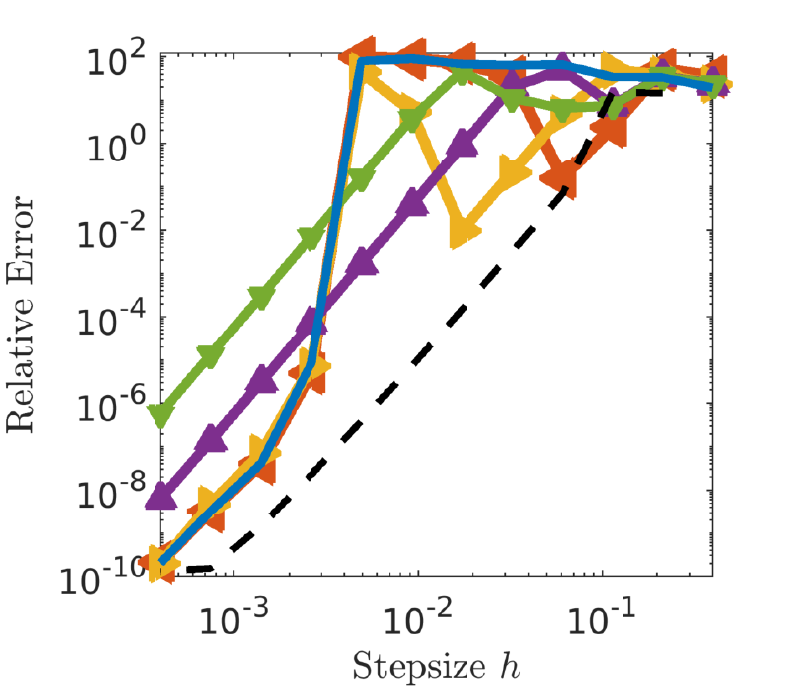} &
			\includegraphics[align=b,width=0.33\linewidth]{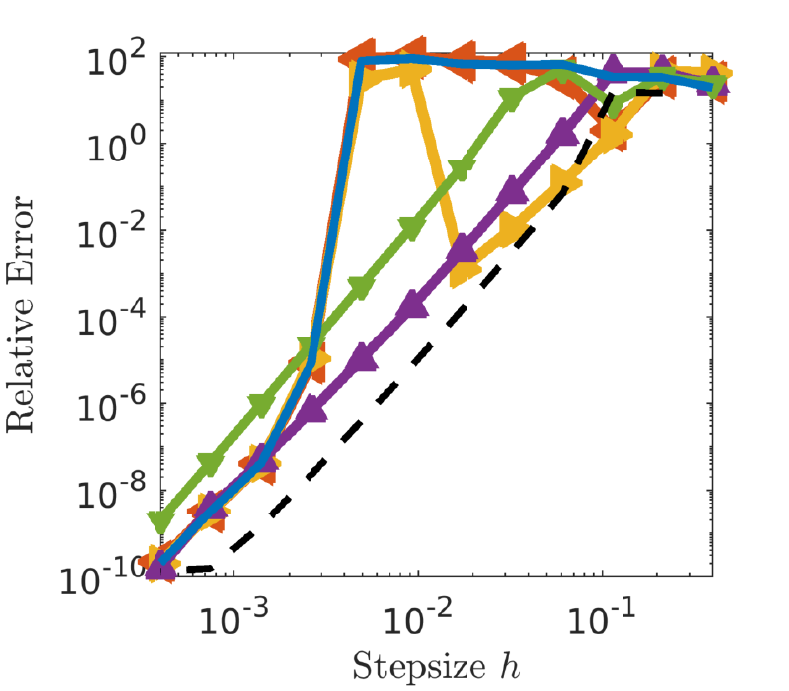} \\
		\end{tabular}
		\end{footnotesize}
		
		\includegraphics[align=b,width=0.9\linewidth,trim={125 190 100 10},clip]{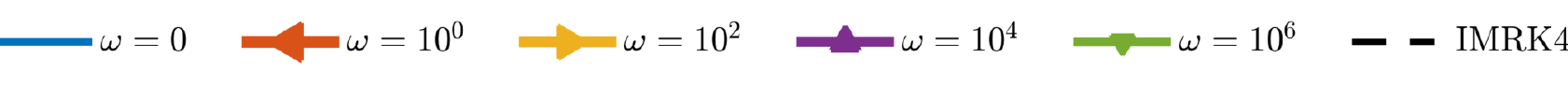}
	
	\end{center}
	\vspace{-1em}
	
	\caption{
		Convergence plots for ERK4 with hyperviscosity added to the linear operator; for convenience we write $\gamma$ in terms of a new parameter $\omega$. Low-order hyperviscosity is unable stabilize the integrator, while a properly selected 8th-order diffusion yields results that are nearly identical to repartitioning.
	}
	\label{fig:erk-hyperviscosity-convergence}

\end{figure}

	To compare the differences between repartitioning and hyperviscosity, we re-solve the ZDS equation using ERK4 with hyperviscosity of orders four, six, and eight. Since ERK4 is a fourth-order method, we take $q=4$. In Figure \ref{fig:erk-hyperviscosity-convergence} we show convergence diagrams for these experiments. We immediately see that hyperviscosity is only effective when a sufficiently high-order spatial derivative is used.  In particular, fourth-order hyperviscosity fails to improve stability for small $\gamma$ and completely damages the accuracy of the integrator for larger $\gamma$. Sixth-order hyperviscosity offers a marginal improvement at coarse stepsizes, but also damages accuracy at fine stepsizes. Eighth-order hyperviscosity with $\gamma=10^{10}$ is the only choice that achieves results comparable to repartitioning.
	
	In summary, repartitioning offers two key advantages. First, it does not require the use of high-order spatial derivatives, and second, it is less sensitive to overdamping. These advantages are both due to the fact that repartitioning does not modify the underlying problem, while  hyperviscosity is only effective if the modified problem (\ref{eq:hyperviscous-semilinear-pde}) closely approximates the original problem (\ref{eq:generic-semilinear-pde}). We discuss both points in more detail below.

	\vspace{1em}
	{\em \noindent Sensitivity to overdamping.} When adding hyperviscosity, it is critical to select the smallest possible $\gamma$ that suppresses instabilities. Selecting larger $\gamma$ causes the solutions of (\ref{eq:generic-semilinear-pde}) and (\ref{eq:hyperviscous-semilinear-pde}) to grow unnecessarily far apart, and leads to a time integration scheme that converges more slowly to the solution of (\ref{eq:hyperviscous-semilinear-pde}). In a convergence diagram, excessive hyperviscosity does not reduce the order-of-accuracy of an integrator, but it will lead to a method with a larger error constant. This phenomenon appears prominently in Figure \ref{fig:erk-hyperviscosity-convergence}, where ERK4 methods with too much hyperviscosity consistently performed worse than all other methods at fine timesteps (e.g. see graphs for $\omega = 10^6$).
	
	 In contrast, second-order and third-order repartitioning are significantly more flexible since they allow for a greater amount of over-partitioning without any significant damage to the accuracy or stability. Excessively large $\epsilon$ values can still cause the stability region separation shown in Figure \ref{fig:exp-overpartitioning}, however such values are unlikely to be used in practice, since they lead to a partitioned linear operator with eigenvalues that have a larger negative real part than imaginary part. Zeroth-order repartitioning is most similar to hyperviscosity since large values of $\epsilon$ also damage the error constant of the method; however, the effects are significantly less pronounced.

	\vspace{1em}
	{\em \noindent Importance of high-order spatial derivatives.} When adding hyperviscosity we must ensure that the small eigenvalues of the modified linear operator $\tilde{\mathbf{L}}$ closely approximate those of $\mathbf{L}$, or we risk altering the dynamics of slow modes. Therefore,  we require a small $\gamma$ for low-order hyperviscous terms. However, this creates a dilemma: choosing a small $\gamma$ may not eliminate the instabilities while choosing a large $\gamma$  damages accuracy. This is exactly why we were not able to efficiently stabilize the ZDS equation using fourth-order hyperviscosity.
			
		In contrast, repartitioning does not require that the small eigenvalues of $\hat{\mathbf{L}}$ closely approximate those of $\mathbf{L}$, since the nonlinear term counteracts any changes. This is perhaps most easily explained by considering the Dahlquist equation (\ref{eq:dispersive-dahlquist}). If $|\lambda| = |\lambda_1 + \lambda_2|$ is small (i.e. the mode is slow), then an exponential integrator will integrate the system accurately so long as $|\lambda_1|$ and $|\lambda_2|$ are also small. Hence, we can freely redistribute the wavenumber between $\lambda_1$ anz $\lambda_2$.  This allows us to repartition using second-order diffusion without loosing accuracy.

\subsection{Long-time integration of Korteweg-de Vries}
\label{subsec:kdv-long-time}

The ZDS equation (\ref{eq:zds}) causes considerable difficulties for unmodified exponential integrators even on short timescales when the dynamics are simple. However, this does not imply that repartitioning is always necessary for solving dispersive equations. In fact, there are many reported results in the literature where exponential integrators converge successfully without requiring any modifications; several examples include \cite{KassamTrefethen05ETDRK4}, \cite{montanelli2016solving}, and \cite{grooms2011IMEXETDCOMP}.%

As discussed in Section \ref{sec:linear_stability}, the instabilities of exponential integrators are very small in magnitude and can therefore go unnoticed for long periods of time. To explore this further, we present a final numerical experiment where we solve the  Korteweg-de Vries (KDV) equation from \cite{zabusky1965interaction}
\begin{align}
	\begin{aligned}
		& \frac{\partial u}{\partial t} = -\left[ \delta \frac{\partial^3 u}{\partial x^3} + \frac{1}{2}\frac{\partial}{\partial x}(u^2) \right] \\ & u(x,t=0) = \cos(\pi x), \hspace{1em} x \in [0,2]
	\end{aligned}
\label{eq:kdv}
\end{align}
where $\delta = 0.022$. The boundary conditions are periodic, and we discretize in space using a 512 point Fourier spectral method. As with the ZDS equation, we solve in Fourier space where the linear operator $\mathbf{L} = \text{diag}(i\delta \mathbf{k}^3)$. 

This exact problem was used in both \cite{buvoli2019esdc} and \cite{buvoli2021epbm} to validate the convergence of ERK, ESDC, and EPBM methods. In the original experiments the equation was integrated  to time $3.6 / \pi$. On these timescales ERK4, ESDC6, and EPBM5 all converge properly and show no signs of instability. %
 We now make the problem more difficult to solve by extending the integration interval to $t=160$. The longer time interval increases the complexity of the solution and allows for instabilities to fully manifest. 

\begin{figure}[h!]

\begin{center}

	\hspace{3em} KDV Solution Plot
	\vspace{1em}
	
	\includegraphics[trim={40 0 55 20},clip,align=b,align=b,width=1\linewidth]{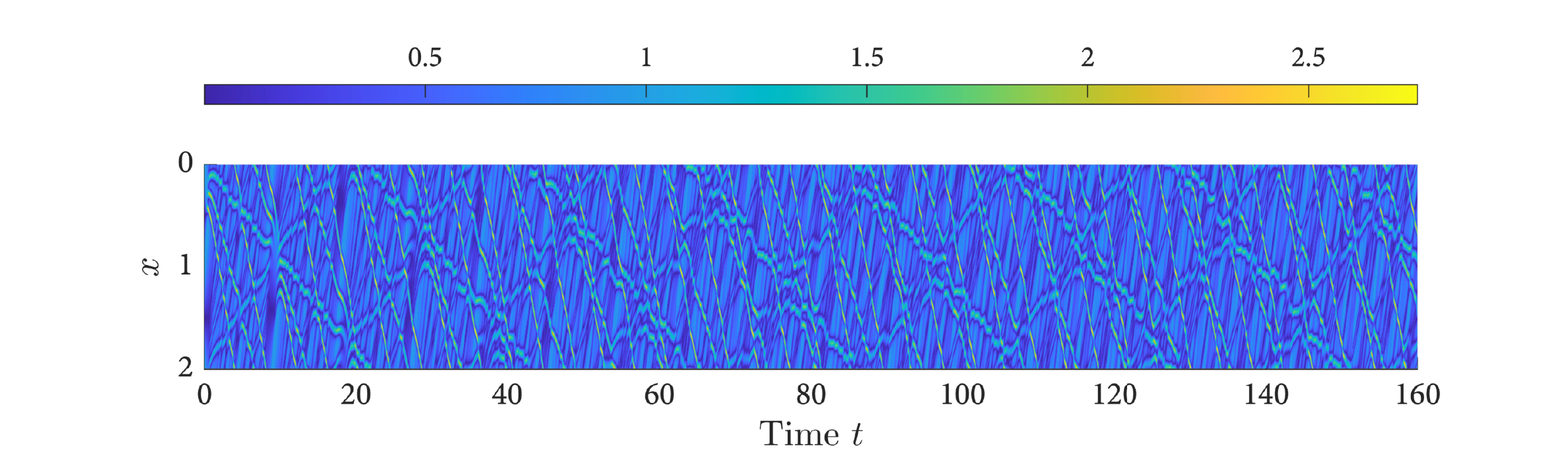}
	
	\vspace{1em}
	\hspace{3em} Error vs. Time
	\vspace{1em}

	\includegraphics[trim={40 0 55 20},clip,align=b,align=b,width=1\linewidth]{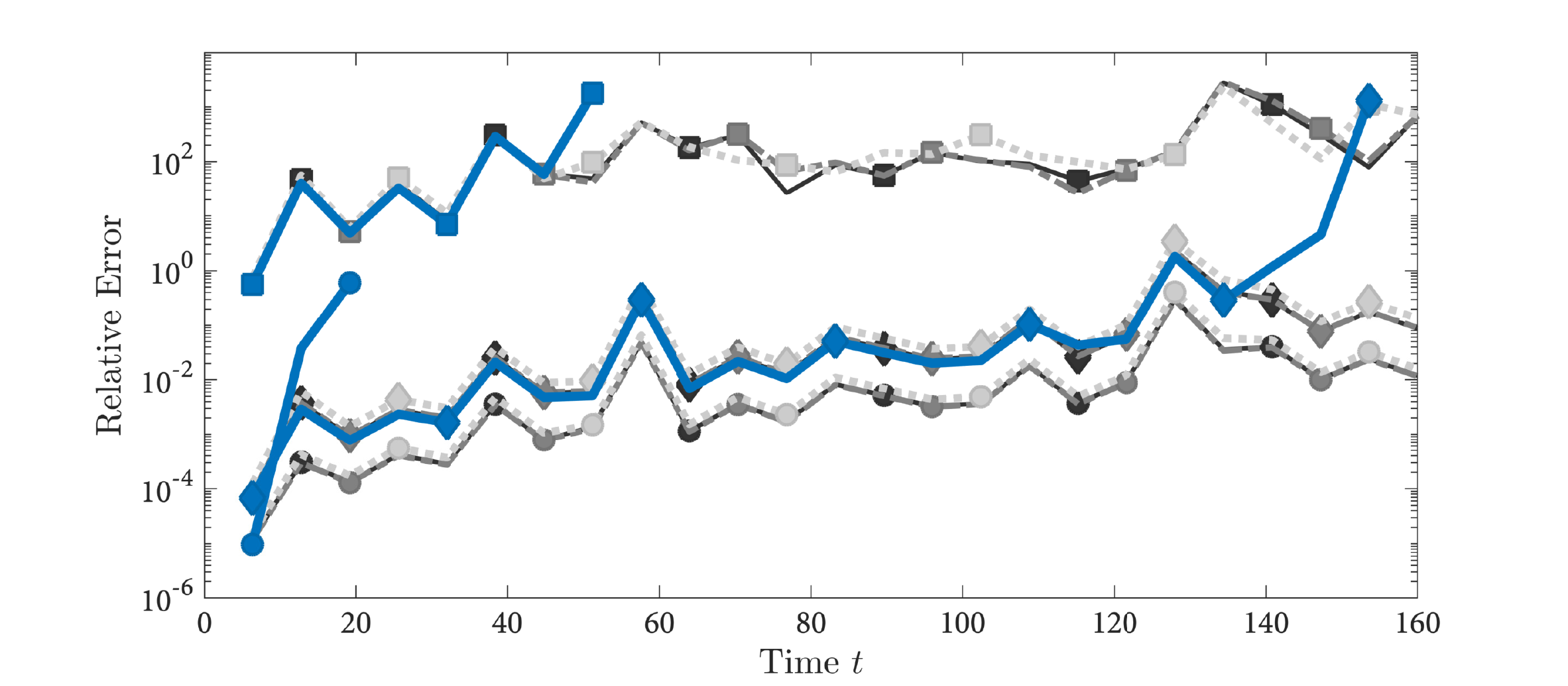}

	\renewcommand{\arraystretch}{2}
	\begin{tabular}[t]{ll} \hline
		Integrator: &
		\begin{tabular}{ccc}
			\includegraphics[trim={250 92 417 7},clip,align=b,height=1em]{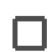} ERK4 &
			 \includegraphics[trim={322 90 344 5},clip,align=b,height=1em]{figures/experiment-kdv/ltc-legend-integrator} ESDC6 &
			 \includegraphics[trim={398 90 270 5},clip,align=b,height=1em]{figures/experiment-kdv/ltc-legend-integrator} ESDC6
		\end{tabular} \\[0em] \hline
		Repartitioning:  &
		\begin{tabular}[t]{ll}
			\includegraphics[trim={90 95 560 8},clip,align=b,align=b,height=1em]{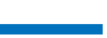} ~ $\gamma = 0$ &
			\includegraphics[trim={165 95 487 8},clip,align=b,height=1em]{figures/experiment-kdv/ltc-legend-repartition.pdf} ~ $\mathbf{D} = -\text{diag}(|\mathbf{k}^3|)$, $\rho = \frac{\pi}{64}$ \\[0em]
				
			\includegraphics[trim={335 95 317 8},clip,align=b,height=1em]{figures/experiment-kdv/ltc-legend-repartition.pdf} ~ $\mathbf{D} = -\text{diag}(\mathbf{k}^2)$, $\rho = \frac{\pi}{3}$ &
			\includegraphics[trim={493 95 159 8},clip,align=b,height=1em]{figures/experiment-kdv/ltc-legend-repartition.pdf} ~ $\mathbf{D} = -\mathbf{I}$, $\gamma = 16$
		\end{tabular} \\[0.5em] \hline
	\end{tabular}
\end{center}
\vspace{-1em}

\caption{Relative error versus time for the KDV equation (\ref{eq:kdv}) solved using ERK4, ESDC6, and EPBM5 with 56000 timesteps. Each integrator is shown using a different line marker. The blue curves denote unpartitioned integrators, while shades of gray denote various repartitionings. The differences in relative error between different repartitionings is very small, so the gray lines almost entirely overlap.}
\label{fig:kdv-long-time-experiment}

\end{figure}

In Figure \ref{fig:kdv-long-time-experiment} we show how the relative error of both unpartitioned and partitoned ERK4, ESDC6, and EPBM5 methods evolves in time. To produce the plot, we run all integrators using 56000 timesteps and compare their outputs to a numerically computed reference solution at 30 equispaced times between $t=0$ and $t=160$. 

On short timescales all unmodified exponential integrators are stable and no repartitioning is required. However, on longer timescales, repartitioning becomes necessary. Moreover, the maximum time horizon before instabilities dominate differs amongst the integrators. The unmodified EPBM5 method is the first integrator to become unstable around $t=20$. The unmodified ERK4 method is more robust and remains stable until approximately time $t=55$, while the unmodified ESDC6 method remains stable across almost the entire time interval. Unlike the ZDS example, the time to instability is now  correlated to the size of the methods stability region.

Adding zeroth-order, second-order, or third-order repartitioning stabilizes all the methods, and does not damage the accuracy in regions where the unmodified integrator converged. Furthermore, the accuracy differences between the three repartitioning strategies is effectively negligible. Lastly, the repartitioning parameters described in the legend of Figure \ref{fig:kdv-long-time-experiment} allow us to compute the solution at even longer times; we tested all methods out to time $t=1000$, after which we did not look further, and found that all partitioned method configurations remained stable.

\section{Summary and conclusions}
\label{sec:conclusion}

In this work we studied the linear stability properties of exponential integrators on stiff equations with no diffusion, and showed that without any modifications, exponential integrators are inherently unstable at large stepsizes. Moreover, we demonstrated that the severity of the instabilities varies across problems and integrator families. %

The main contribution of this paper is the introduction of a simple repartitioning approach that resolves the stability issues and improves convergence for stiff non-diffusive equations. Compared to the common practice of adding hypervisocity, repartitioning provides better accuracy and does not require introducing high-order spatial derivatives. Furthermore, there is no need to solve a modified equation; instead repartitioning introduces a pre-specified amount of numerical dissipation directly into an exponential integrator. 

The applications of this work are two-fold. First, repartitioning is a useful tool for any practitioner experiencing stability issues when applying exponential integrators on hyperbolic or dispersive equations. The simplest way to apply repartitioning is to solve an equation using a range of different repartitioning terms until the desired result is achieved. Alternatively, one can study the linear stability properties of an integrator, and then compute the spectrum of the equation's linear operator to determine a suitable repartitioning term and constant $\epsilon$.

The second application of this work relates to exponential integrators and parallel-in-time frameworks where effect of instabilities becomes much more significant. In followup works we will apply this repartitioning strategy within the Parareal and PFASST frameworks to solve non-diffusive equations.

\section*{Acknowledgements}
The work of Buvoli was funded by the National Science Foundation, Computational Mathematics Program DMS-2012875. The work of Minion was supported by the U.S. Department of Energy, Office of Science, Office of Advanced Scientific Computing Research, Applied Mathematics program under contract number DE-AC02005CH11231.

\bibliographystyle{siam}
\bibliography{references_nourl,pint}

\appendix
\section{Method coefficients}
\label{app:coefficients}

\begin{itemize}[leftmargin=*]
		
\item { \bf ERK4} is a fourth-order exponential Runge-Kutta method given by
\begin{align*}
	K_0 &= N(t_n, y_n) \\
	K_1 &= N(t_n + h/2, \varphi_0(h\mathbf{L}/2)y_n + \tfrac{1}{2} \varphi_1(h/2 \mathbf{L}) K_0 ) \\
	K_2 &= N(t_n + h/2, \varphi_0(h\mathbf{L}/2)y_n + \left( \tfrac{1}{2} \varphi_1(h/2 \mathbf{L}) - \varphi_2(h/2 \mathbf{L}) \right) K_0 + \varphi_2(h/2 \mathbf{L})K_1 )
\\
	K_3 &= N(t_n + h, \varphi_0(h\mathbf{L})y_n + \left( \varphi_1(h \mathbf{L}) - 2 \varphi_2(h \mathbf{L}) \right) K_0 + 2 \varphi_2(h \mathbf{L}) K_2)\\
	y_{n+1} &= \varphi_0(h\mathbf{L})y_n + 	(\varphi_1(h\mathbf{L}) - 3 \varphi_2(h\mathbf{L}) + 4 \varphi_3(h\mathbf{L})) K_1 \\
	& ~~~ + (2 \varphi_2(h\mathbf{L}) - 4 \varphi_3(hL)) (K_2 + K_2) +  (- \varphi_2(hL) - 4 \varphi_3(h\mathbf{L})) K_3
\end{align*}

\item { \bf ESDC6}  is a sixth-order exponential spectral deferred correction method. To write down the formula for ESDC6, we first define the quadrature nodes
\begin{align*}
	t_{n,j} = t_n + h \eta_j, && j=1,\ldots, 4
\end{align*}
for $\eta_j = \{ 0, \tfrac{1}{2} - \frac{\sqrt{5}}{10}, \tfrac{1}{2} + \frac{\sqrt{5}}{10}, 1\}$ and let $h_{n,j} = t_{n,j+1} - t_{n,j}$. The pseudocode for ESDC6 is:

\vspace{1em}
\renewcommand{\arraystretch}{1.25}
\hspace{1em} \begin{tabular}{l}
        $Y^{[k]}_{n,1} = y_n$, ~$N^{[k]}_{n,j} = N(t_{n,j}, Y^{[k]}_{n,j})$ \\
        \text{{\bf for} j = 1 to $3$} \\        
        \hspace{1.5em} $Y^{[1]}_{n, j+1} = \varphi_0(h_{n,j}\mathbf{L}) Y^{[1]}_{n, j} + h_{n,j} \varphi_1(h_{n,j} \mathbf{L}) N_{n,j+\gamma}^{[1]}$ \\
        \text{{\bf for} k = 1 to $6$}\\
        \hspace{1.5em} \text{{\bf for} j = 1 to $3$}\\
        \hspace{3em}	$Y^{[k+1]}_{n,j+1} = \varphi_0(h_{n,j} \mathbf{L}) Y^{[k+1]}_{n,j} + h_{n,j}\varphi_1(h_{n,j} \mathbf{L}) \left[ N_{n,j+\gamma}^{[k+1]} - N_{n,j+\gamma}^{[k]} \right] + I_{n,j}^{[k]}$ \\[0.5em]
        $y_{n+1} = Y^{[m+1]}_{n,p}$ \\[1em]     
\end{tabular}
To write the exponential integral terms $I_{j,k}$ we first define
\begin{align*}
   	\tau_{j,i} = \frac{t_{n,i} - t_{n,j}}{h_{n,j}}, \quad i,j = 1, \ldots, 4, \\
   \mathbf{V}(j)_{c,d} = (\tau_{j,c})^{d-1}, \quad 	\mathbf{V}(j) \in \mathbf{R}^{4\times 4},  
\end{align*}
then,
\begin{align}
    I_{n,j}^{[k]} &= h_{n,j} \sum_{\nu = 1}^p  \varphi_{\nu}(h_{n,j} \mathbf{L}) \mathbf{b}^{[k]}_{\nu} \quad \text{where} \quad \mathbf{b}^{[k]}_{\nu} = \sum_{l=1}^p a^{(\nu)}_{j,l} \hspace{0.1em} N_{n,l}^{[k]}.
    \label{eq:exponential_integral_a}
\end{align}
where the weights $a_{j,l}^{(\nu)} = \nu! \mathbf{V}(j)^{-1}_{\nu + 1, l}$.

\item {\bf EPBM5} is composite method based on an exponential polynomial block method that accepts five inputs ${y_j^{[n]} \approx y(t_n + r z_j)}$ and produces five outputs ${y_j^{[n+1]} \approx y(t_n + r z_j  + r\alpha)}$ where the nodes $\{z_j\}_{j=1}^5$ are $\{ -1, -\eta^{+}, -\eta^{-}, \eta^{-}, \eta^{+} \}$ for $\eta^{\pm} = \sqrt{\tfrac{3}{7} \pm \tfrac{2}{7} \sqrt{\tfrac{6}{5}}}$.
The output computation is
\begin{align}
	y^{[n+1]}_j = \varphi_0(r \eta_j(\alpha) \mathbf{L}) y^{[n]}_1 + r \sum_{k=1}^{4} \eta_j^{k} \varphi_{k}(r \eta_j(\alpha) \mathbf{L}) \mathbf{v}_k, && j = 1, \ldots 5,
	\label{eq:epbm-timestep-computation}
\end{align}
where $\eta_j(\alpha) = z_j + \alpha + 1$ and the vectors $\mathbf{v}_j$ are given by 
\begin{align*}
	\begin{bmatrix}
		\mathbf{v}_1 \\
		\mathbf{v}_2 \\
		\mathbf{v}_3 \\
		\mathbf{v}_4 \\
	\end{bmatrix}
	=
	\begin{bmatrix}
		(w_1^+ + u_1^+) 	& (-w_1^- + u_1^-) 	& (w_1^- + u_1^-) 	& (-w_1^+ + u_1^+) \\
		(-w_2^- - u_2) 	& (w_2^+ + u_2)		& (-w_2^+ + u_2)		& (w^-_2 - u_2) \\
		(w_3^- + u_2)	& (-w_3^+ - u_2)		& (w_3^+ - u_2)		& (-w_3^ +- u_2) \\
		-w_4^+			& w_4^+				& w_4^+				& w_4^-
	\end{bmatrix}
	\begin{bmatrix}
		\mathbf{N}_2 \\
		\mathbf{N}_3 \\
		\mathbf{N}_4 \\
		\mathbf{N}_5 \\
	\end{bmatrix}
\end{align*}
where $\mathbf{N}_j = N(t_n + rz_j, y^{[n]}_j)$ and the constants $w_j^{\pm}$ and $u_j^{\pm}$ are: 
\begin{align*}
	w_1^{\pm} &= \tfrac{\sqrt{75 \pm 4\sqrt{30}}}{12}	 & w_2^{\pm} &= \tfrac{\sqrt{10170 \pm 1104 \sqrt{30}}}{24} \\
	w_3^{\pm} &= \tfrac{7\sqrt{1350 \pm 180 \sqrt{30}}}{24}	& 	w_4^{\pm} &= \tfrac{7 \sqrt{150 \pm 20\sqrt{30}} }{8}  \\
	u^{\pm}_1 &= \tfrac{3 \pm \sqrt{30}}{12} & u_2 &= \tfrac{7 \sqrt{30}}{24}
\end{align*}
EPBM5 first advances the solution using $\alpha=1$, and then corrects the new solution using $\alpha=0$. If we denote the right-hand-side of (\ref{eq:epbm-timestep-computation}) as $M_j(\alpha, t_n, y^{[n]})$, then the composite method EPBM5 can be written as
\begin{align*}
	\begin{aligned}
		\tilde{y}_j^{[n+1]} &= M_j(1, t_n, y^{[n]}) \\
		y_j^{[n+1]} &= M_j(0, t_n + r, \tilde{y}^{[n+1]})	
	\end{aligned}
	&& j = 1, \ldots, 5.
\end{align*}
The required nonlinear function evaluations $\mathbf{N}_j = N(t_n + rz_j, y^{[n]}_j)$ and $\mathbf{N}_j = N(t_n + rz_j + r, \tilde{y}^{[n]}_j)$ can each be evaluated in parallel.

\end{itemize}

\end{document}